\documentclass[a4paper,11pt]{article}
\usepackage{amsmath}
\usepackage{amsfonts}
\usepackage{titlesec}
\usepackage{amssymb}
\usepackage{amsthm}
\usepackage{enumitem}
\usepackage[frozencache=true,cachedir={./}]{minted} 
\usepackage{geometry}
\usepackage{color}

\usepackage[utf8]{inputenc}

\usepackage{mathtools}
\usepackage{algorithm}
\usepackage{array}
\makeatletter
\newcommand{\thickhline}{%
	\noalign {\ifnum 0=`}\fi \hrule height 1pt
	\futurelet \reserved@a \@xhline
}
\newcolumntype{"}{@{\hskip\tabcolsep\vrule width 1pt\hskip\tabcolsep}}
\makeatother

\usepackage{listings}

\newcommand{\ci}{\mathrm{i}} 

\usepackage{framed}
\usepackage{graphicx}
\usepackage{caption, subcaption}
\usepackage{hyperref}
\usepackage{float}

\newtheorem{theorem}{Theorem}[section]

\newtheorem{remark}[theorem]{Remark}

\newcommand{\quotes}[1]{``#1''}
\newcommand{\uOD}{u^{\text{\tiny LOD}}}

\newcommand{\VLOD}{V_{H,\text{\tiny LOD}}^{h,\ell}}
\newcommand{\VOD}{V_H^{\text{\tiny LOD}}}
\newcommand{\LOD}{_{\text{\tiny LOD}}}
\newcommand{\R}{\mathbb{R}}

\newcommand{\D}{\mathcal{D}}

\newcommand{\Iptr}{\text{I}_\text{ptr}}
\newcommand{\J}{\text{J}_{\boldsymbol{\omega}}}
\newcommand{\K}{\text{K}_{\boldsymbol{\omega}}}

\newcommand{\sub}{{\fontfamily{lmtt}\selectfont sub}}

\newcommand{\dx}{\hspace{2pt}\mbox{d}x}
\newcommand{\Corr}{\mathcal{Q}}

\newcommand{\lambdags}{\lambda_\text{gs}}
\newcommand{\lambdagsLOD}{\lambda_{\text{gs}}^{\LOD}} 
\newcommand{\ugs}{u_\text{gs}}
\newcommand{\ugsLOD}{u_{\text{gs}}^{\LOD}} 
\newcommand{\Vgs}{V}
\newcommand{\Egs}{E}
\newcommand{\betags}{\beta}

\newcommand{\ac}{a_{\text{c}}}
\newcommand{\aV}{a_V}

\newcommand{\PLODLtwo}{P^{\LOD}_{L^2}} 
\newcommand{\PHLOD}{P^{\LOD}} 
\newcommand{\Ltwo}[2]{(#1,#2)}

\definecolor{myBlue}{RGB}{113,104,238} 
\definecolor{myGreen}{RGB}{154,205,50} 
\definecolor{myGreen2}{RGB}{114,175,30} 
\definecolor{myRed}{RGB}{180,50,50}  
\definecolor{myOrange}{RGB}{225,92,22} 
\definecolor{lgray}{RGB}{200,200,200} 
\definecolor{llgray}{RGB}{155,155,155} 
\definecolor{myForestGreen}{RGB}{34,139,34}
\definecolor{lila}{rgb}{0.5,0.0,0.5}

\definecolor{mycolor1}{rgb}{0.00000,0.44700,0.74100}%

\usepackage{diagbox}

\usepackage{bbm} 

\begin{document}

\begin{center}
{\LARGE 
A two level approach for simulating Bose--Einstein condensates by localized orthogonal decomposition}
\end{center}

\begin{center}
{\large Christian D\"oding\footnote[1]{Institute for Numerical Simulation, University Bonn, DE-53115 Bonn, Germany, \\ e-mail: \textcolor{blue}{doeding@ins.uni-bonn.de}.}, 
Patrick Henning\footnote[2]{Department of Mathematics, Ruhr University Bochum, DE-44801 Bochum, Germany, \\ e-mail: \textcolor{blue}{patrick.henning@rub.de}.} and 
Johan W\"arneg{\aa}rd\footnote[3]{Department of Applied Physics and Applied Mathematics, Columbia University, New York, NY 10027, U.S., \\ e-mail: \textcolor{blue}{jpw2176@columbia.edu}.}}\\[2em]
\end{center}

\begin{center}
{\large{May 23, 2024}}
\end{center}

\begin{abstract}
In this work, we consider the numerical computation of ground states and dynamics of single-component Bose--Einstein condensates (BECs). The corresponding models are spatially discretized with a multiscale finite element approach known as Localized Orthogonal Decomposition (LOD). Despite the outstanding approximation properties of such a discretization in the context of BECs, taking full advantage of it without creating severe computational bottlenecks can be tricky. In this paper, we therefore present two fully-discrete numerical approaches that are formulated in such a way that they take special account of the structure of the LOD spaces. One approach is devoted to the computation of ground states and another one for the computation of dynamics. A central focus of this paper is also the discussion of implementation aspects that are very important for the practical realization of the methods. In particular, we discuss the use of suitable data structures that keep the memory costs economical. The paper concludes with various numerical experiments in $1d$, $2d$ and $3d$ that investigate convergence rates and approximation properties of the methods and which demonstrate their performance and computational efficiency, also in comparison to spectral and standard finite element approaches.
\end{abstract}
%
\vspace{12pt}

\vspace{12pt}

\vspace{12pt}
\noindent

\section{Introduction}
In this paper we consider the numerical approximation of the stationary and the time-dependent Gross--Piteaevskii equation by advanced finite element techniques. The Gross--Pitaevskii equation has applications in quantum physics where it is used to model so-called Bose--Einstein condensates (BECs) which form when a dilute bosonic gas of particles is cooled down to near the absolute zero temperature \cite{PiS03}. In this extreme state of matter all particles within the condensate are in the same quantum state and are indistinguishable from each other so that the BEC itself behaves like a single quantum particle. This state, sometimes called a \quotes{super-atom}, gives access to quantum mechanical phenomena on an observable scale such as superfluidity and justifies why BECs became a highly active research area in modern quantum physics; see e.g. \cite{CornellWiemanNobel,Nature,KetterleNobel,CornellWieman99,Ketterle01} and the references therein. An underlying mathematical model that describes the formation of BECs was first derived by Gross \cite{Gro61} and Pitaevskii \cite{Pit61} and later mathematically justified by Lieb et al. \cite{LSY00,LSY01}.

We are interested in numerical simulations of real-world models in the context of BECs. For that reason it is required to consider two problems of different mathematical flavor. The first one is the time-dependent Gross--Pitaevskii equation (GPE),
\begin{align}
\label{preliminary:GPE}
\ci \partial_t u = -\tfrac{1}{2} \triangle u + V u + \beta |u|^2u ,
\end{align}
which describes the time evolution of the quantum mechanical wave function $u$ of the BEC. The function $V = V(x)$ models a magnetic or optical potential that traps the condensate in the physical experiment, and the constant $\beta$ is a repulsion parameter so that the nonlinearity $\beta|u|^2u$ models the interaction between the particles within the condensate. The second problem we have to tackle is the nonlinear Gross--Pitaevskii eigenvalue problem (GPEVP)
\begin{align}
\label{preliminary:GPEVP}
	\lambda u = -\tfrac{1}{2}\triangle u + Vu + \beta |u|^2 u,
\end{align}
which characterizes energy extrema of the wave function. The unique minimal solution $u$, minimal in terms of $\lambda$, is called the ground state whereas all other solutions are called excited states. Besides their strong mathematical connection, there is also high interest in studying both problems in view of the application. The reason for this is the standard procedure followed in experiments with BECs; c.f. \cite{BaoNum,Nature}. In a first step, a BEC is created by cooling down the gas of bosons to near absolute zero temperature. During this process, the condensate is held in place by a magnetic or optical trapping potential. As a result of this first stage, the condensate is prepared in the ground state w.r.t. the current trapping potential. In a second stage, the potential is tweaked or removed and the subsequent dynamics of the condensate are studied. In this stage the wave function of the condensate is described by the solution of the time-dependent GPE w.r.t. the tweaked trapping potential. \\
The main goal of this work is to present an advanced computational approach for simulating both stages of the aforementioned physical process and therefore for tackling both mathematical problems. The approach is based on a generalized finite element approximation in space known as the  \emph{Localized Orthogonal Decomposition} (LOD) which is a multiscale technique that was first introduced by M{\aa}lqvist and Peterseim \cite{HeM14} in the context of highly heterogeneous diffusion problems and was further successfully applied to the GPEVP \cite{HeMaPe14,HeP21} and the GPE \cite{Superconv}. In this work we further extend the in  \cite{HeMaPe14,HeP21,Superconv} achieved results and design two novel methods that combine the LOD with suitable iterative solvers depending on the nature of the considered problem, i.e., the GPEVP or the time-dependent GPE. The designed fully-discrete methods come with several advantages in view of its approximation properties due to the LOD discretization. However, realizing them in an efficient way is a nontrivial task since the methods require the assembly of nonlinear terms in each iteration which may harm the theoretical advantages of the methodology in view of the computational costs. Therefore, we introduce a modification of the desired methods by adding an additional projection on the LOD space into the method. As we shall then demonstrate in numerical experiments this modification significantly decreases the computational coast and does not impose further substantial errors so that we achieve a very good overall performance of these novel methods. \\
To solve the nonlinear eigenvalue problem GPEVP, which is equivalent to the minimization of the associated energy functional, an adaptive Riemannian optimization method is used as an interative solver and that is based on changing the metric in each iteration. The numerical algorithm of the iteration was originally proposed and analyzed in \cite{SobolevGradient} and enjoys global convergence properties, which makes it very attractive in practice to obtain reliable approximations. For the time-dependent GPE, a continuous Galerkin (cG) ansatz is used for the time discretization, first introduced in combination with classical Lagrange finite elements in space by Karakashian and Makridakis \cite{KaM99} for the cubic nonlinear Schr\"odinger equation and further analyzed in \cite{DoHen22}. The resulting time-stepping method is energy-conserving 
and 
it is of convergence order $\mathcal{O}(\tau^{2q})$ with respect to the time-step size $\tau$ and where $q$ denotes the polynomial degree of the time ansatz functions (cf. \cite{KaM99,DoHen22}). We observe numerically the same order of convergence when we combine the time-integrator with the more problem-specific LOD space. \\
The algebraic structure of the equation \eqref{preliminary:GPE} comes with several time invariants that convey information about the physical state of the system on a global scale. With this in mind, it comes as no surprise that in order to achieve reliable numerical methods it is of utmost importance to conserve and accurately represent these quantities in the discrete setting. We recall that the accuracy with which these quantities are reproduced in the discrete setting is subject to a suitable spatial discretization whereas their conservation is subject to a suitable time discretization. As for the time discretization, the choice of conservative time-integrator schemes over non-conservative ones can have a tremendous advantage in terms of accuracy. This observation has been confirmed in various numerical experiments (cf. \cite{NLSComparison,Sanz-Serna}).  Among others, mass conservative time discretizations have been studied in \cite{Wang,zouraris_2001}, time integrators that are mass conservative and symplectic are investigated in  \cite{Akrivis1991,HeM17,Tou91,Verwer1984}, energy conservative time discretizations in \cite{KaM99} and time discretizations that preserve mass and energy simultaneously are addressed in \cite{Akrivis1991,BaC12,Besse,BDD18,CCW20,NonlinearCN,H1Est,Sanz-SernaNLCN}. \\
As for the spatial discretization, a popular choice is to use a spectral approach which, for time integrations, often relies on splitting methods. These have been analyzed in for example \cite{BaoNum,Lub08,Tha12b}. Other common choices include finite differences, finite element methods or recently holistic discretizations \cite{BunderRoberts}. Recent cross comparisons between different approaches have been made in \cite{ANTOINE20132621, BaoNum,Spectral,NLSComparison}. In particular, in high regularity regimes with smooth potentials and regular domains, the exponential convergence of a spectral approach leads to methods that are computationally very efficient. On the downside, they typically perform poorly in low regularity regimes such as in the case of a low regularity potential, which shall also be exemplarily illustrated in the present paper. Here, the method of LOD presents a hitherto unique opportunity to bridge this gap as it achieves high spatial accuracy even for discontinuous potentials, as well being computationally highly competitive for smooth solutions.

In many physically relevant settings it can also be reasonable to include the error in the initial discrete energy as an additional measure for the accuracy of a method. An example for such a setting was illustrated in \cite{Superconv}, where excrescent discrete energy in the system was converted into artificial kinetic energy that lead to spurious solutions. One way of achieving high accuracy with respect to 
the initial energy at low computational cost was addressed in \cite{Superconv} by using the aforementioned LOD. The LOD allows to build a discrete approximation space that is adapted to the underlying differential equation by incorporating problem-specific information (explained in section \ref{section-LOD}). With this technique it is possible to approximate the initial value to 6th order accuracy with respect to the mesh size in terms of the time invariants of the equation. This superconvergence property only requires moderate regularity of the initial value $u_0$, which is compatible with physically relevant problems where $u_0$ minimizes a certain energy. This observation makes the LOD extremely attractive in the context of the Gross--Pitaevskii equation and its conserved quantities, not only for the time-dependent equation as in \cite{Superconv}, but also for the Gross--Pitaevskii eigenvalue problem \cite{HeMaPe14,HeP21}. For further reading on the LOD we refer to the textbook by M{\aa}lqvist and Peterseim \cite{BookLODAxelDaniel} and the survey article \cite{ActaNumericaLOD}. Basic implementation aspects are discussed in \cite{ENGWER2019123} and for more recent developments on the LOD we exemplarily refer to \cite{HaP23,HeKeMa20,HMW21,KMP22,LMP21,Maier21,MaV22,WuZh22}.

We finish the introduction with a brief outline of the paper. In section \ref{section:math-setup} we set up the mathematical framework we are working in and formulate the aforementioned two problems in the context of BECs. Section \ref{section-LOD} is then devoted to the spatial discretization by Localized Orthogonal Decomposition. There we start with the basic construction of the LOD space, apply it to the GPE and then discuss the approximation properties regarding ground states and time invariants. The section is closed with a description of the practical implementation of the LOD discretization. In section \ref{section:CompEVP} we then derive the fully-discrete method that we propose for solving the GPEVP using a damped inverse iteration based on LOD discretization. In section \ref{section:high-order-cG-FEM} we state the fully-discrete method for solving the time-dependent GPE using the energy conserving time integrator. We finish the paper in section \ref{section:numerical-experiments} with numerical results presenting the performance and efficiency of our methods in several experiments. Further notes on the efficient implementation of the numerical methods is given in the appendix in section \ref{section:implementation-notes}. 

\section{Mathematical setting of the Gross--Pitaevskii equation}
\label{section:math-setup}
In this section we describe the mathematical framework of the Gross--Pitaevskii equation and  formulate the two problems for the formation of ground states and for the dynamic evolution of BECs. Throughout this paper, we shall denote by ${\mathcal{D}} \subset \R^d$ (for $d=1,2,3$) a bounded convex domain with polyhedral boundary. The complex Hilbert space $L^2(\D) = L^2(\D,\mathbb{C})$ of complex-valued square-integrable functions is equipped with the scalar product $\Ltwo{ v}{ w }=\int_{\mathcal{D}} v(x) \, \overline{w(x)} \dx$ where $\overline{w}$ is the complex-conjugate of $w$ and we denote by $\| \cdot \|$ the induced norm for brevity. As usual, $H^1_0(\D) = H^1_0(\D,\mathbb{C})$ denotes the Sobolev space of complex-valued, weakly differentiable functions with zero trace on $\partial \mathcal{D}$ and $L^2$-integrable partial derivatives. In the derivation of the problems \eqref{preliminary:GPE} and \eqref{preliminary:GPEVP} for modeling BECs, the energy of the system is a quantity that plays an important role. It is given by the Gross-Pitaevskii energy functional
\begin{align} \label{GPenergy}
	E(u) = \int_{\D} \tfrac{1}{2}|\nabla u|^2 + V|u|^2 + \tfrac{\beta}{2} |u|^4 \dx
\end{align}
for $u \in H^1_0(\D)$ where we assume that the trapping potential is given by a bounded, real-valued, and non-negative function, i.e., $V \in L^\infty(\D,\R_{\ge 0})$ and that the repulsion parameter $\beta \in \R$ is non-negative. The Gross-Pitaevskii energy functional \eqref{GPenergy} is derived from the many body problem of a cold and weakly interacting bosonic gas by mean-field approximation and a Hatree-ansatz for the wave function of the condensate. We refer to \cite{BaoNum,Dynamics,Gro61,LSY01,Pit61,PiS03} for precise details on the physical derivation and validation of the Gross--Pitaevskii energy functional.

\subsection{The Gross--Pitaevskii eigenvalue problem and ground states}

We now derive the Gross--Pitaevskii eigenvalue problem \eqref{preliminary:GPEVP} from the energy functional \eqref{GPenergy} and specify the notion of ground states of BECs. In particular, we are interested in finding the state of lowest energy of a BEC (in the given trapping potential $V$), i.e., the state that minimizes the Gross--Pitaevskii energy \eqref{GPenergy} under the normalization constraint $\|u \| = 1$. Hence, we seek for $\ugs \in H^1_0(\D)$ such that
\begin{align} \label{definition-groundstate}
\ugs =\underset{u \in \mathbb{S}}{\mbox{arg\hspace{2pt}min}}\hspace{2pt} \Egs(u)
\qquad \mbox{where } \mathbb{S}:= \{  u \in H^1_0(\D) \mbox{ and } \| u\| =1 \}.
\end{align}
We call $\ugs$ the \emph{ground state} of the BEC (in the given configuration) and it can be shown that the problem is well-posed and $\ugs$ is unique up to sign \cite{CCM10}. It is possible to characterize $\ugs$ equivalently through a nonlinear eigenvalue problem. This perspective can be both beneficial for the analysis, as well as for numerical computations. To reformulate \eqref{definition-groundstate} as an eigenvalue problem, we simply consider the corresponding Euler-Lagrange equations of the form $\Egs^{\prime}(\ugs) = \lambdags \, S^{\prime}(\ugs)$, where $S(u):=\| u \|^2 - 1$ represents the constraint functional (as we demand $S(\ugs)=0$) and $\lambdags$ is the Lagrange multiplier for the constraint. Computing the Fr\'echet derivatives of the two functionals we recover the GPEVP in variational form as 
\begin{align*}
\lambdags  \Ltwo{ \ugs }{ v } = \tfrac{1}{2}\Ltwo{ \nabla \ugs }{ \nabla v } + \Ltwo{ \Vgs \ugs}{v} + \betags \Ltwo{|\ugs|^2\ugs }{ v } \qquad \mbox{for all } v\in H^1_0(\D),
\end{align*}
with the constraint $\ugs \in \mathbb{S}$. In this setting, the eigenfunction $\ugs$ associated with the smallest eigenvalue $\lambdags$ is indeed the ground state in the sense of \eqref{definition-groundstate}. This nontrivial observation was first made in \cite{CCM10}. Consequently, we can write the GPEVP compactly as: find the smallest eigenvalue $\lambdags>0$ and corresponding eigenfunction $\ugs \in \mathbb{S}$ such that in the sense of distributions on $H^1_0(\D)$:
\begin{align}
\label{definition-groundstate-2}
	\lambdags \ugs =  -\tfrac{1}{2}\triangle \ugs + \Vgs \ugs+ \betags |\ugs|^2 \ugs.
\end{align}
We note that the ground state $\ugs$ is (up to sign) identical to the solution of \eqref{definition-groundstate}. Natural regularity statements about $\ugs$ that are important for optimal convergence rates of our numerical approximations are given in \cite{HeP21}. 
\begin{remark}
There are various iterative approaches for finding minimizers of \eqref{definition-groundstate}, such as methods based on Riemannian optimization \cite{APS22,AnD14,ALT17,BaD04,DaK10,DaP17,HSW21}, based on self consistent field iterations (SCF) \cite{DiC07,UJR21}, or the \quotes{$J$-method} \cite{AHP21b,JKM14}. 
In this paper, we apply the damped inverse iterations proposed in \cite{SobolevGradient} due to its fast and global convergence properties when the damping parameter is chosen adaptively. Details of the method are given later in section \ref{section:CompEVP}.
\end{remark}

\subsection{The time-dependent Gross--Pitaevskii equation and its invariants}
\label{subsection:time-invariants}

Next, let us turn to the time-dependent GPE describing the dynamics of BECs on a relevant time interval $[0,T]$. We assume that the BEC is at time $t = 0$ in a state which is described by a wave function $u_0 \in H^1_0(\D)$. The time evolution of the BEC is then determined by the Schr\"odinger equation $\ci \partial_t u = \tfrac{1}{2} E'(u)$ with the Gross--Pitaevskii energy functional $E$ from \eqref{GPenergy}. Taking the Fr\'echet derivative of the energy functional into account leads to the time-dependent GPE formulated in variational form as
\begin{align} \label{weakGPE}
	\Ltwo {\ci \partial_t u}{ v} = \tfrac{1}{2}\Ltwo{ \nabla u }{ \nabla v } + \Ltwo{ V u}{v} + \beta \Ltwo{|u|^2u }{ v } \qquad \mbox{for all } v\in H^1_0(\D).
\end{align}
Hence we can formulate the time-dependent Gross-Pitaevskii problem as follows: find 
$$
u \in L^{\infty}([0,T],H^1_0({\mathcal{D}})) 
\qquad
\mbox{and}
\qquad
\partial_t  u \in L^{\infty}([0,T],H^{-1}({\mathcal{D}}))
$$
such that $u(0)=u_0 \in H^1_0(\D)$ and such that for almost every $t \in (0,T]$
\begin{eqnarray}
\label{model-problem}\ci \partial_t u = -\tfrac{1}{2} \triangle u 
+  V \hspace{1pt} u + \beta |u|^2 \hspace{1pt} u.
\end{eqnarray}
Note that such a solution automatically belongs to $C([0,T],L^2(\D))$ and therefore the notion of the initial value is well-defined. Moreover, the initial value $u_0 \in H^1_0(\D)$ can be chosen arbitrarily in the general setting. However, as motivated by our physical application, we shall mainly consider the physically reasonable case that $u_0=\ugs$ is a ground state in a modified configuration in the sense of problem \eqref{definition-groundstate-2}, i.e., $\ugs$ from \eqref{definition-groundstate-2} is the ground state trapped in a different potential than used for the time-dependent problem \eqref{model-problem}. This being said, it is worth to mention that the maximal existence time of solutions to the GPE \eqref{model-problem} depends crucially on the initial value (in particular its regularity and the size of $u_0$ in the $H^1$-norm), so that unphysical initial values easily lead to a chaotic behavior after short times that ultimately leads to the nonexistence of solutions. All aspects of existence and uniqueness of solutions to the GPE \eqref{model-problem} are nicely elaborated in the textbook by Cazenave \cite{Cazenave}. In particular, the GPE is at least locally well-posed, i.e., for $T$ sufficiently small. When talking about solutions and discussing convergence results for numerical methods, we shall always assume that problem \eqref{model-problem} admits a \emph{unique solution} that is \emph{sufficiently regular}.

Since the problem \eqref{model-problem} is derived as a Hamiltonian system from the Gross--Pitaevskii energy functional \eqref{GPenergy} its solution $u$ is automatically energy-preserving over time, that is
\begin{align*}
	E(u(t)) = E(u_0) \quad \text{for all } t \in [0,T].
\end{align*}
For example, this can be derived directly from the variational form \eqref{weakGPE} when it is tested with $v = \partial_t u$ and imaginary parts are taken. In contrast to that, when \eqref{weakGPE} is tested with $v = u$, one infers the second fundamental conservation law
\begin{align*}
	\| u(t) \|^2 = \| u_0 \|^2 \quad \text{for all } t \in [0,T]
\end{align*}
which is the preservation of mass. Other conservation laws hold true in less general settings such as the conservation of momentum and conservation of the mass center in the absence of the trapping potential. For details and more possible conservation laws we refer examplarily to \cite{BaoNum,Superconv}. From the numerical point of view it can be extremely important to choose numerical methods that approximate these quantities with high accuracy in order to capture the correct behavior and dynamics of BECs in numerical experiments; see e.g. \cite{HeM17,NLSComparison,Sanz-Serna}.

\section{Spatial discretization by Localized Orthogonal Decomposition}
\label{section-LOD}
As sketched in the introduction, our spatial discretization is based on the Localized Orthogonal Decomposition. In the following we will briefly describe the construction of LOD spaces and why it can be beneficial to use them in various applications. After that we review the approximation properties of these spaces in the context of the GPEVP and the time-dependent GPE and we describe one of the central computational issues in this context, namely the efficient assembly of nonlinear terms.

Before diving into the explicit construction of the spaces, we shall introduce some preliminaries that hold throughout the whole manuscript.
The first essential ingredient for constructing a LOD space is a \quotes{coarse} finite element mesh $\mathcal{T}_H$ that corresponds to a shape-regular and quasi-uniform partition of the computational domain $\D$ with mesh size $H>0$. In the following, we shall restrict the type of considered partitions to being simplicial (i.e., intervals in 1D, triangles in 2D and tetrahedra in 3D), though the construction also directly generalizes to quadrilateral partitions. By a \emph{\quotes{coarse} mesh} we mean that the number of elements in $\mathcal{T}_H$ is small compared to what would be usually necessary to obtain reasonable approximations in corresponding conventional finite element spaces. The aspect of coarseness will become clearer when describing the construction in more detail below.\\

Next, we consider a standard P1 finite element space on the mesh $\mathcal{T}_H$ given by
$$
V_H := \{ v \in C(\overline{\D}) \cap H^1_0(\D) | \hspace{3pt} v_{\vert K} \in \mathbb{P}^1(K) \mbox{ for all } K \in \mathcal{T}_H \}.
$$
We shall later modify this space (while keeping its dimension fixed), by correcting it based on problem-specific information entering through an elliptic differential operator. A crucial component of the construction is the choice of a set of quantities of interest which are represented by pairwise linear independent functionals $q_j : H^1_0(\D) \rightarrow \R$. Given such a set, the basic task is to construct \quotes{cheap} numerical approximations $u_H^{\LOD}$ that essentially conserve all the quantities of interest for the exact solution $u$, i.e., with the property $q_j(u_H^{\LOD}) = q_j(u)$ for all $q_j$. For linear problems, exact conservation is often possible, whereas for nonlinear problems (like the GPE), we have to content with almost conservation up to high accuracy. In the context of the GPE, quantities of interest that were identified as suitable are {\it local weighted averages} based on the coarse space $V_H$, cf. \cite{HeMaPe14,HeP21,Superconv}. More precisely, we select
$$
q_j(v) := ( \varphi_j , v )_{L^2(\D_j)},
$$
where $\varphi_j \in V_H$ denotes the (conventional) nodal basis function associated with the (interior) mesh node $z_j$ and which is supported on a patch $\D_j \subset \mathcal{D}$. Before stating how this is exploited in the construction of a space, we require one more notation, which is the so-called detail space $W\subset H^1_0(\D)$ that contains functions that do not change the quantities of interest when added to $H^1$-functions and in particular to elements of $V_H$. It is given by
\begin{align*}
W := \{ w \in H^1_0(\D) | \hspace{3pt} q_j(w)= 0 \mbox{ for all } j \,\}.
\end{align*}
Note that $W$ is just the kernel of the $L^2$-projection onto the finite element space $V_H$ as its definition implies $\Ltwo{ v_H }{ w  } =0 $ for all $v_H \in V_H$ and all $w \in W$.

\subsection{The ideal Localized Orthogonal Decomposition}
\label{subsection:idealized_LOD}
To explain how to correct the finite element space $V_H$ based on an elliptic differential operator and the quantities of interest encoded through the space $W$, we shall study an idealized setting. For that let $\mathcal{L}$ be a symmetric elliptic differential operator represented by a (symmetric) bilinear form $a(\cdot,\cdot)$ on $H^1_0(\D)$ and let $f \in L^2(\D)$ be a given source term. With this, consider the equation seeking $u\in H^1_0(\D)$ with
\begin{align}
\label{ref-problem}
a(u,v) = \Ltwo{ f}{ v}  \qquad \mbox{for all } v\in H^1_0(\D).
\end{align}
A suitable discrete solution space is now constructed based on the $a(\cdot,\cdot)$-orthogonal projection on $W$ denoted by $\Corr : H^1_0(\D) \rightarrow W$, i.e., for $v\in H^1_0(\D)$, the image $\Corr v \in W$ is given as the solution to
\begin{align}
\label{idealized-def-C}
a( \Corr v , w) = a(v,w) \qquad \mbox{for all } w\in W.
\end{align}
The operator $\Corr$ is also called the \emph{correction operator}, because we use it to \quotes{correct} $V_H$ in the sense that the ideal approximation space is defined by
\begin{align}
\VOD := \{ v_H - \Corr v_H | \hspace{3pt} v_H \in V_H\},\label{LODdef}
\end{align}
or compactly $\VOD = (I - \Corr) V_H$. Obviously, $\VOD$ has the same dimension as the coarse space $V_H$. Furthermore, Galerkin approximations in $\VOD$ conserve the quantities of interest. To see this, let $\uOD \in \VOD$ denote the solution to
\begin{align}
\label{LOD-problem}
a( \uOD , v) = \Ltwo{ f }{ v } \qquad \mbox{for all } v \in \VOD.
\end{align}
Testing in \eqref{ref-problem} with $v \in \VOD$ and subtracting it from \eqref{LOD-problem} yields the Galerkin orthogonality
\begin{align}
\label{idealize-Galerkin-orth}
a( \uOD - u , v) = 0 \qquad \mbox{for all } v \in \VOD.
\end{align}
As $\VOD = (I - \Corr) V_H$ is the $a(\cdot,\cdot)$-orthogonal complement of $W$ (which directly follows from \eqref{idealized-def-C}), we conclude from \eqref{idealize-Galerkin-orth} that $ \uOD - u \in W$ and consequently
\begin{align*}
 q_j( \uOD - u ) = 0, \quad \mbox{or equivalently } \quad  q_j( \uOD ) = q_j( u ),
\end{align*}
for all quantities of interest $q_j$. By the choice of $q_j$ as $q_j(v)=\Ltwo{ \varphi_j }{ v }$, we also obtain error estimates that only depend on the regularity of the source term $f$. To be precise, let $s=0,1,2$, where $s=0$ if $f\in L^2(\D)$ and $s=1,2$ if $f\in H^1_0(\D) \cap H^s(\D)$, then it holds
\begin{align*}
\| \uOD - u \| + H \, \| \uOD - u \| _{H^1(\D)} \le C\, H^{s+2} \| f \|_{H^s(\D)},
\end{align*}
where the constant $C$ does not depend on $u$ or its regularity. For a detailed derivation of the estimate we refer to \cite{Superconv}. The essential observation to make here is that for elliptic PDEs it is possible to construct a problem-specific low-dimensional space $\VOD$ that admits super-approximation properties that do not depend on the regularity of the solution itself or possible structural variations of $u$. In particular, the underlying mesh $\mathcal{T}_H$ is generic and only needs to compensate $\| f\|_{H^s(\D)}$, but it does not need to resolve the variations of $u$. The latter aspect is indirectly taken care of through the space $\VOD$.

\subsection{Application to the Gross--Pitaevskii equation} 
\label{subsection:superconvergence-results}
In order to approximate the solutions of the Gross--Pitaevskii problems \eqref{definition-groundstate-2} and \eqref{model-problem} using the LOD discretization in space (as described in the previous subsection) we have to specify the bilinear form $a(\cdot,\cdot)$ on $H^1_0(\D)$. For this purpose, we neglect the nonlinear term in the GPE and may choose the bilinear form $a(\cdot,\cdot)$ in \eqref{ref-problem} as
\begin{align}
\label{choice-LOD-a}
\aV(u,v):=\tfrac{1}{2}\Ltwo{ \nabla u }{ \nabla v } + \Ltwo{ V u }{ v}
\end{align}
which is induced by the main part of the GPE and, in particular, incorporates the trapping potential $V$. Alternatively, it is also possible to include only parts of the potential $V$ or to fully drop it which results in the choice
\begin{align}
\label{choice-LOD-a0}
\ac(u,v):=\Ltwo{ \nabla u }{ \nabla v }
\end{align}
to which we refer as the canonical bilinear form. There are a few practical differences regarding the performance of the LOD based on the choices \eqref{choice-LOD-a} and \eqref{choice-LOD-a0}. First of all, if $V$ is a rough potential with low regularity, it is typically reasonable to include it in the construction of $\VOD$ to obtain optimal convergence rates, cf. \cite{HeP21}. If $V$ is smooth, the guaranteed convergence rates are the same whether $V$ is included or not, even though the absolute errors might differ quite significantly. The major practical difference between \eqref{choice-LOD-a} and \eqref{choice-LOD-a0} lies in the computation of the corresponding LOD basis functions. If $V$ is periodic or not included at all, then the LOD basis functions $(I - \Corr)\varphi_j$ have to be computed only for a few nodes $z_j$, whereas the remaining basis functions are essentially obtained by translations. On the contrary, if $V$ is heterogenous, then it is necessary to compute each LOD basis separately by solving an elliptic problem (cf. \eqref{idealized-def-C}), which causes considerable costs in the initial phase (\quotes{offline phase}). Last but not least, if $V$ is included, then $(I - \Corr)\varphi_j$ has often a more complex structure. This calls for more accurate quadrature rules when computing products of (gradients of) LOD basis functions in mass and stiffness matrices.

Many of these aspects will become clearer after describing the method in full detail together with its implementation. We shall also revisit these issues in the numerical experiments in section \ref{section:numerical-experiments} to highlight them again and to support them by CPU times.\\[0.6em]
In the following we shall briefly review the relevant approximation properties of LOD spaces in the context of the Gross--Pitaevskii equation.

\subsubsection{Ideal LOD-approximation of ground states}
\label{subsubsection-LOD-ground-states}
We consider the GPEVP for the ground state eigenfunction $\ugs$ and we recall that is is equivalent to the energy minimization problem \eqref{definition-groundstate}. 
In order to approximate $\ugs$, let $\VOD$ be the LOD space constructed as described in section \ref{subsection:idealized_LOD} based on a coarse P1 finite element space $V_H$ and on the bilinear form $a(\cdot,\cdot)$ which is either chosen as in \eqref{choice-LOD-a} or as in \eqref{choice-LOD-a0}. A numerical approximation of $\ugs$ is then given by a discrete minimizer in the LOD space, i.e.,
\begin{align*}
\ugsLOD  =\underset{v \in \mathbb{S} \cap \VOD}{\mbox{arg\hspace{2pt}min}}\hspace{2pt} \Egs(v).
\end{align*}
As $\ugsLOD$ and $\ugs$ can be only unique up to sign, we can select them such that they are consistent in the sense that
\begin{align*}
\Ltwo{ \ugs }{ \ugsLOD } \ge 0.
\end{align*}
Since we assumed $\Vgs\in L^{\infty}(\D)$ to be non-negative and $\betags \ge 0$, we have the following a priori error estimates for the LOD approximations that were proved in \cite{HeP21}:
\begin{align*}
\| \ugs - \ugsLOD \| + H \| \ugs - \ugsLOD \|_{H^1(\D)} \le C H^4.
\end{align*}
Here the constant $C>0$ depends on the potential $\Vgs$, the ground state eigenvalue $\lambdags$, the constant $\betags$, the domain $\mathcal{D}$, the maximum value of the ground state $\| \ugs \|_{L^{\infty}(\D)}$ and the mesh regularity of $\mathcal{T}_H$. For the corresponding ground state energy one has the estimate
\begin{align*}
\left| \Egs( \ugs ) - \Egs( \ugsLOD) \,\right| \le C H^6
\end{align*}
and for the ground state eigenvalue
\begin{align*}
\left|  \lambdags - \lambdagsLOD \right| \le C H^6 \qquad \mbox{where }  \lambdagsLOD := \int_{\D} \tfrac{1}{2} |\nabla \ugsLOD|^2 + \Vgs | \ugsLOD |^2 + \betags | \ugsLOD|^4 \dx.
\end{align*}
The above result shows that in the LOD space, we obtain superconvergence in various errors. For comparison, in the standard P1 finite element space $V_H$ (which has the same dimension as $\VOD$), the $H^1$-error only converges with the order $\mathcal{O}(H)$ and the $L^2$-error with the order $\mathcal{O}(H^2)$. This is the same for the errors in the ground state energy and in the ground state eigenvalue which also both converge with the order $\mathcal{O}(H^2)$ (compared to $\mathcal{O}(H^6)$ in $\VOD$). The superior approximation properties of the LOD space make it very attractive for solving the Gross--Pitaevskii eigenvalue problem. If $\ugsLOD$ is selected as an initial value in an energy-conserving discretization of the time-dependent equation \eqref{model-problem}, then we can guarantee that the energy of the discrete approximations accurately resembles the exact energy on arbitrarily large time scales.

\subsubsection{Ideal LOD-approximation of time invariants}
\label{subsubsection-LOD-time-invariants}
Let us now turn towards the question how accurate time invariants of the GPE, such as the energy and mass, can be approximated in the LOD space. This question particularly arises when we assume that the initial value $u_0$ in \eqref{model-problem} is given by an analytical function in $H^1_0(\D)$. Hence, since we aim to solve \eqref{model-problem} numerically using the LOD discretization in space we first need to project the initial value onto the LOD space $\VOD$. But the solutions of the GPE may be very sensitive to perturbations in energy (or mass) and therefore this error w.r.t. energy in the initial value caused by the spatial discretization is of high relevance. \\
In general, let us assume that a sufficiently smooth function $u \in H^1_0(\D)$ is given and that the space $\VOD$ is constructed based on the bilinear form $a(\cdot,\cdot)$, selected either as in \eqref{choice-LOD-a} or as in \eqref{choice-LOD-a0}. With this, a suitable approximation of $u$ in the LOD space is given by its $a(\cdot,\cdot)$-orthogonal projection, i.e., $\PHLOD (u) \in \VOD$ is the solution to
\begin{align*}
a( \PHLOD(u) ,  v) = a( u , v) \qquad \mbox{for all } v \in \VOD.
\end{align*}
The central question is now how accurately $\PHLOD(u)$ approximates the important time invariants introduced in section \ref{subsection:time-invariants} and if we can obtain similar rates as for the ground state. A positive answer was given in \cite{Superconv}, where it was proved that for sufficiently smooth $u$ and for a sufficiently smooth potential $V$, the energy (cf. \eqref{GPenergy}) and the mass are both approximated with 6th order accuracy, i.e.,
\begin{align}
\label{optimalrates-energy-mass-LOD}
\big| E(u) - E(\PHLOD(u)) \big| \le C H^6 \qquad \mbox{and} \qquad \big| \| u \|^2 - \| \PHLOD(u) \|^2 \big| \le C H^6,
\end{align}
where the constant $C$ depends on $u$, $V$, $\beta$ and parameters that characterize the quasi-uniformity and the regularity of the mesh. Similar results are obtained for the momentum and the center of mass in the absence of the trapping potential; see \cite{Superconv}. We note however that the optimal rates \eqref{optimalrates-energy-mass-LOD} can be only formally guaranteed for the solution $u$ of the time-dependent GPE if the potential fulfills $V\in H^2(\D)$, though this can be relaxed to some extend for suitable initial values (cf. Remark \ref{remark-wellprepared-gs} below).

In conclusion we can summarize that, both for analytically given (sufficiently regular) initial values $u_0$ and precomputed ground states, the LOD spaces allow to approximate the energy with 6th order accuracy. Preserving it over time is then subject to a suitable time discretization of the GPE. 

\begin{remark}[Well-prepared initial values for optimal accuracy]
\label{remark-wellprepared-gs}
As mentioned before, the initial values for the time-dependent GPE are often selected as ground states. Now let us consider such a setting and denote by $V^{\text{gs}}$ the potential used for describing the ground states and by $V^{\text{t}}$ the potential in the time-dependent GPE. Let us assume that both potentials split into a regular and into a rough contribution, i.e., $V^{i}=V_0^{i}+V_d^{i}$ with $V_0^i \in H^2(\D)$ and $V_d^{i} \in L^{\infty}(\D)$ where $i= \text{gs} , \text{t}$. If both problems are discretized with the same LOD space, then optimal convergence orders for the energy can only be guaranteed if the rough contributions of the potentials are identical (i.e. $V_d^{\text{gs}}=V_d^{ \text{t}}$) and if $V_d^{\text{gs}}$ is included in the construction of the LOD space. In this case, we can interpret $\ugs$ as a well-prepared initial value for the time-dependent GPE because $V_d^{\text{t}}$ will only mildly reduce the time-regularity for $u(t)$ (cf. \cite[Appendix A]{NonlinearCN}). A rigoros error analysis in this setting is however still open.
\end{remark}

\subsection{Practical approximation of the ideal LOD-spaces}
Before we can start discussing implementation details, we need to revisit the idealized LOD spaces and how they are approximated in practice by a two level discretization. We will do this exemplarily for the space constructed based on the bilinear form $a(\cdot,\cdot) = \aV(\cdot,\cdot)$ from \eqref{choice-LOD-a} but we can easily replace it by \eqref{choice-LOD-a0}. \\
An approximation of the ideal space $\VOD$ is necessary since the canonical set of basis functions of $\VOD$ consists of functions $(I - \Corr) \varphi_j$, where $\varphi_j$ are the nodal basis functions of $V_H$ and $\Corr: H^1_0(\D) \rightarrow W$ denotes the correction operator given by \eqref{idealized-def-C}. As the computation of $\Corr \varphi_j$ requires the solution of an elliptic problem in the infinite-dimensional Hilbert space $W$ (cf. \eqref{idealized-def-C}), we need to discretize it on a second finite element mesh $\mathcal{T}_h$. The mesh $\mathcal{T}_h$ needs to be sufficiently fine and is usually selected as a subdivision of the coarse mesh $\mathcal{T}_H$ in order to simplify the implementation. The corresponding P1 finite element space $V_h$ on $\mathcal{T}_h$ thus contains the coarse space, i.e., $V_H\subset V_h\subset H^1_0(\D)$. With this, we could solve the corrector problem \eqref{idealized-def-C} in the discrete detail space
$$
W_h := V_h \cap W,
$$
leading to a discrete corrector problem of the form
\begin{align}
\label{ideal-discrete-corrector-problem}
\mbox{find } \Corr_h \varphi_j \in W_h: \qquad 
a(  \Corr_h \varphi_j , w_h  ) = a( \varphi_j , w_h  ) \quad \mbox{for all } w_h \in W_h.
\end{align}
However, in this case a second issue remains, that is, that even though the right hand side of \eqref{ideal-discrete-corrector-problem} is local (due to the local support of $\varphi_j$), the solution $\Corr_h \varphi_j $ has (in general) support on the whole domain $\D$, which would make the repeated solving for each basis function $\varphi_j$ very expensive. Luckily, due to the particular solution space $W_h$, the functions $\Corr_h \varphi_j $ are exponentially decaying (in units of the coarse mesh size $H$) outside of the support of $\varphi_j$, i.e., outside the nodal patch $\D_j$. This property was first proved in \cite{LocalElliptic} and justifies that the global corrector problems can be truncated to small local patches. This makes the solving feasible and the expenses for assembling an approximation of $\VOD$ remain economical. From the view point of approximation theory it is important that the local basis functions are constructed such that they still form a partition of unity in order to avoid numerical pollution \cite{MeB96}.

In the following we shall present a suitable localization strategy to approximate $\Corr_h \varphi_j$. The strategy is commonly associated with the classical LOD and was suggested in \cite{HeM14, HeP13}. For that, let $v_H \in V_H$ be an arbitrary function in the coarse space $V_H$ (in practice selected as a nodal basis function $\varphi_j$) and let $K \in \mathcal{T}_H$ be a coarse simplex. For the localization we require a truncation parameter $\ell \in \mathbb{N}_{>0}$ that characterizes the size of local patches. To be precise, we define the $\ell$-layer patch around $K$ iteratively by 
\begin{align*}
S_\ell(K) & : = \bigcup \{ \hat{K} \in \mathcal{T}_H | \hspace{3pt}\hat{K} \cap S_{\ell-1}(K) \not= \emptyset \} \qquad \mbox{and} \qquad S_0(K) := K.
\end{align*} 
This means that $S_\ell(K)$ consists of $K$ and $\ell$ layers of grid elements around it. This also implies that $S_\ell(K)$ has a diameter of order $\mathcal{O}(\ell H)$. The restriction of $W_h$ to $S_\ell(K)$ is given by
$$
W_h(S_\ell(K)) := W_h \cap H^1_0(S_\ell(K))
$$
and the restriction of the bilinear form $a(\cdot,\cdot)$ to some local subset $S \subset \mathcal{D}$ by
\begin{align*}
a_{S}(u,v) := \tfrac{1}{2} ( \nabla u , \nabla v )_{L^2(S)} +  ( V u ,  v )_{L^2(S)}.
\end{align*}
With all these notations, the $K$-localized corrector $\Corr_{h,K}^{\ell} : V_H \rightarrow W_h(S_\ell(K))$ is given by
\begin{align}
\label{def-local-truncated-corrector}
a_{S_\ell(K)} ( \Corr_{h,K}^{\ell}v_H , w_h  ) = a_K( v_H , w_h) \qquad \mbox{for all } w_h \in  W_h(S_\ell(K)).
\end{align}
Note that computing $\Corr_{h,K}^{\ell}v_H$ involves solving a local finite element problem on the (small) domain $S_\ell(K)$. 
The global approximation of the idealized corrector $\Corr_h$ (cf. \eqref{ideal-discrete-corrector-problem}) is then defined by
\begin{align*}
\Corr_{h}^{\ell}v_H := \sum_{K \in \mathcal{T}_H} \Corr_{h,K}^{\ell}v_H.
\end{align*} 
It is easy to see that if $\ell$ is large enough such that $S_\ell(K)=\mathcal{D}$, then we have
\begin{align*}
\Corr_{h}^{\ell}v_H = \Corr_{h} v_H,
\end{align*} 
i.e., the approximated corrector coincides with the ideal corrector. In general, the error between $\Corr_{h}^{\ell}$ and $\Corr_{h}$ can be quantified by the aforementioned exponential decay as
\begin{align}
\label{corr-estimate-ell}
\| \Corr_{h}^{\ell}v_H - \Corr_{h} v_H \|_{H^1(\D)} \le C \,  \exp(-\rho \ell ) \, \| v_H \|_{H^1(\D)},
\end{align}
where $C>0$ and $\rho>0$ are generic constants that depend on the mesh regularity and on the coercivity and continuity constants of $a(\cdot,\cdot)$. For a proof we refer to \cite{HeM14,HeP13}. \\
\noindent
In consequence, we can define the approximated LOD space by
\begin{align*}
\VLOD := \{ (I - \Corr_{h}^{\ell})v_H | \hspace{3pt} v_H \in V_H \},
\end{align*}
which is spanned by the basis functions $(I - \Corr_{h}^{\ell}) \varphi_j$. Naturally, the localization introduced through $\Corr_{h}^{\ell}$ becomes a source of numerical error. However, as indicated by \eqref{corr-estimate-ell}, the approximation properties of $\VLOD$ are the same as for $\VOD$, up to an exponential error term of the form $ \exp(-\rho \ell )$. For rigorous proofs of this claim, we refer to \cite{BookLODAxelDaniel}. As a practical remark, selecting $\ell$ proportional to $|\log(H)|$ with a suitable prefactor allows to fully reproduce the high convergence orders stated in sections \ref{subsubsection-LOD-ground-states} and \ref{subsubsection-LOD-time-invariants}.

Hence, in practice we can work with the space $\VLOD$ without sacrificing significant accuracy. The advantages of $\VLOD$ are apparent: On the one hand, the basis functions $(I - \Corr_{h}^{\ell}) \varphi_j$ have a local support (of order $|\log(H)|H$) which is important to obtain sparse system matrices after the spatial discretization, both for the eigenvalue problem and for the time-dependent problem. On the other hand, a basis function  $(I - \Corr_{h}^{\ell}) \varphi_j$ is cheap to compute. It requires to solve for $\Corr_{h,K}^{\ell} \varphi_j$ for every coarse element $K$ that is in the nodal patch $\D_j$. Each of the problems \eqref{def-local-truncated-corrector} is small (only formulated on $S_\ell(K)$) and they are all independent from each other and can be solved in parallel. Also note that problem \eqref{def-local-truncated-corrector}
 can be written as a saddle point problem where the constraint entering through the space $W_h(S_\ell(K))$ is represented by a Lagrange multiplier. The solution space for the saddle point problem now becomes the unconstraint space $V_h(S_\ell(K))$. This is exploited in practical computations.
 
The implementation of the above constructed LOD space is described comprehensively in \cite{ENGWER2019123} which is why we will not repeat it here. There is only one aspect that we would like to mention as it becomes important when discussing the implementation of the numerical schemes for the Gross--Pitaevskii equation: Once the LOD basis functions are computed, they are stored and processed through a matrix. To be a bit more specific, let $\varphi_i$, $i=1,\dots,N_H$ denote the (nodal) basis functions of $V_H$, let $\psi_j$, $j=1,\dots,N_h$ denote the (nodal) basis functions of $V_h$, and let $\mathbf{P} \in  \R^{N_H\times N_h}$ denote the interpolation matrix containing the values of each coarse basis function $\varphi_i$, in every fine node, so that
\begin{align*}
	\varphi_i = \sum_{j = 1}^{N_h} \mathbf{P}_{ij} \psi_j \qquad \mbox{for all } 1 \le i \le N_H.
\end{align*}
Furthermore, we can store the coefficients of the functions $\Corr_{h}^{\ell} \varphi_i \in V_h$ (for each $1\le i \le N_H$) row-wise in a corrector matrix $\mathbf{Q} \in  \R^{N_H\times N_h}$ such that
\begin{align*}
	\Corr_{h}^{\ell} \varphi_i = \sum_{j = 1}^{N_h} \mathbf{Q}_{ij} \psi_j \qquad \mbox{for all } 1 \le i \le N_H.
\end{align*}
Then the representation of the LOD-basis functions on the fine mesh is given by the \emph{LOD-matrix} $\boldsymbol{\Phi}^{\LOD} = \mathbf{P} - \mathbf{Q}$, such that
\begin{align*}
	\varphi^{\LOD}_i = \sum_{j = 1}^{N_h} \boldsymbol{\Phi}^{\LOD}_{ij} \psi_j \qquad \mbox{for all } 1 \le i \le N_H
\end{align*}
is a basis of the LOD space $\VLOD$. With this, the standard stiffness matrix can be assembled easily through matrix multiplications, e.g. for the stiffness matrix in the LOD space $\mathbf{A}^{\LOD}_{ij} = \Ltwo{\nabla \varphi^{\LOD}_i}{\nabla \varphi^{\LOD}_j}$ we have
$\mathbf{A}^{\LOD} = (\mathbf{P} - \mathbf{Q}) \mathbf{A} (\mathbf{P} - \mathbf{Q})^{\top}$
where $\mathbf{A}_{ij} = \Ltwo{\nabla \psi_i}{\nabla \psi_j}$ is the stiffness matrix on the fine mesh. Similarly, we can assemble the mass matrix and the matrix representing the potential term in the GPE. However, this simple strategy fails once nonlinear terms are considered as they will appear in the GPE. In the next subsection we present a practical solution for this issue.

\subsection{Treatment of nonlinear terms}
\label{subsection:treatment-nonlinearity-with-LOD}
The super-approximation properties of the LOD-space reviewed in section \ref{subsection:superconvergence-results} can be further leveraged to substantially speed up the nonlinear computations. This is however not obvious since the GPE involves a cubic nonlinearity of the form $\Ltwo{ |u|^2 u}{ v }$ which has to be discretized in the LOD space. With a standard discretization (and suitable linearization), this would essentially involve a product of the form $\Ltwo{ |v^{\LOD}|^2 \varphi_i^{\LOD}}{ \varphi_m^{\LOD}}$, where $v^{\LOD} \in \VLOD$ is a function in the LOD space and $\varphi_i^{\LOD} = (I - \Corr_{h}^{\ell}) \varphi_i$ again denotes the LOD basis function belonging to the node with index $i$. The necessary precomputations for easily evaluating such a term would require to assemble and store a 4-valence tensor in the LOD space, which easily exceeds the computing and memory resources due to the curse of dimensionality.

To this end, we introduce the \emph{sparse} 3-valence tensor 
\begin{align} \label{def_omega}
	\boldsymbol{\omega}_{ijk} := \Ltwo{ \varphi_i^{\LOD}\varphi_j^{\LOD}}{ \varphi_k^{\LOD}}, \quad 1 \le i,j,k \le N_H
\end{align} 
and replace throughout expressions of the type $\Ltwo{ |v|^2u}{w}$ in the GPE with $\Ltwo{ \PLODLtwo(|v|^2)u}{w }$, where $\PLODLtwo$ is the $L^2$-projection onto the LOD space, i.e., $\PLODLtwo(|v|^2) \in \VLOD$ solves
\begin{align}
\label{PLODL2}
\Ltwo{ \PLODLtwo(|v|^2) }{ \varphi_i^{\LOD} } = \Ltwo{ |v|^2 }{ \varphi_i^{\LOD} } \qquad \mbox{for all } 1 \le i \le N_H.
\end{align}
With this notation, the projection $\PLODLtwo(|v\LOD|^2)$ for an LOD function $v\LOD$ is easily computed using $\boldsymbol{\omega}$. In fact, the LOD-mass matrix describing the left hand side of \eqref{PLODL2} is available (typically preassembled as sketched in the previous subsection) and the entries of the right hand side vector $\Ltwo{ |v\LOD|^2 }{ \varphi_i^{\LOD}}$ are obtained from $\boldsymbol{\omega}$. To be precise, if $v\LOD= \sum_k \boldsymbol{\alpha}_k \,\varphi_k^{\LOD}$ for $\boldsymbol{\alpha}_k \in \mathbb{C}$, then the entries of the right hand side vector in \eqref{PLODL2} can be expressed as
\begin{align*}
\Ltwo{ |v\LOD|^2 }{ \varphi_i^{\LOD} } & = \Ltwo{ \sum_k\sum_j \boldsymbol{\alpha}_k \overline{\boldsymbol{\alpha}_j} \,\, \varphi_k^{\LOD} \varphi_j^{\LOD} }{ \varphi_i^{\LOD} } 
= \sum_{k\leq j }
\Re(\boldsymbol{\alpha}_k \overline{\boldsymbol{\alpha}_j}) \hspace{2pt} (2-\delta_{kj}) \hspace{2pt}
\Ltwo{  \varphi_k^{\LOD} \varphi_j^{\LOD} }{ \varphi_i^{\LOD} } \\
& =  \sum_{k\leq j} \Re(\boldsymbol{\alpha}_k \overline{\boldsymbol{\alpha}_j}) \hspace{2pt} (2-\delta_{kj}) \hspace{2pt} \boldsymbol{\boldsymbol{\omega}}_{kji}, \qquad\mbox{where } \boldsymbol{\omega}_{kji}= \Ltwo{ \varphi_k^{\LOD} \varphi_j^{\LOD} }{ \varphi_i^{\LOD} }.
\end{align*}
Now that \eqref{PLODL2} is algebraically assembled, we can cheaply solve the arising system of equations to obtain a vector representation of $\PLODLtwo(|v\LOD|^2)$, say,
\begin{align*}
\PLODLtwo(|v\LOD|^2) =  \sum_k \boldsymbol{\rho}_k \,\varphi_k^{\LOD} \qquad
\mbox{for some } \boldsymbol{\rho}\in\mathbb{C}^{N_H}.
\end{align*}
With this, nonlinear terms of the form $\Ltwo{ \PLODLtwo(|v\LOD|^2) u\LOD }{ \varphi_i^{\LOD} }$ are then straightforwardly assembled as
\begin{align*}
\Ltwo{ \PLODLtwo(|v\LOD|^2) u\LOD }{ \varphi_i^{\LOD} } 
=  \sum_{kj} \boldsymbol{\rho}_k  \boldsymbol{\alpha}_j \, \boldsymbol{\omega}_{kji}.
\end{align*}
We note that the threefold symmetry of $\boldsymbol{\omega}$ implies that only nonzero elements for which $i\leq j\leq k$, must be computed. Using the 3-valence tensor instead of the canonical 4-valence tensor reduced the preassembly time and the memory demand tremendously and ultimately makes the implementation of the LOD-based numerical methods efficient. In the Appendix \ref{section:implementation-notes} we give more details on the implementation of the 3-valence tensor $\omega$.

Remarkably, the error introduced through replacing $|v|^2$ by $\PLODLtwo(|v|^2)$ at all occurrences seems to be only of order $\mathcal{O}(H^8)$ and is hence negligible. Heuristically, this can be motivated by the following simple calculation. First, we observe that for any $u,v\in H^1_0(\D)$ it holds:
\begin{align*}
\int_\D |u|^2 v -\PLODLtwo(|u|^2)v \dx = \int_{\D} (|u|^2-\PLODLtwo(|u|^2))(v-\PLODLtwo(v)) \dx.
\end{align*} 
If $u,v$ and $V$ are sufficiently smooth, that is $(-\tfrac{1}{2}\Delta   + V) |u|^2  \in H^1_0(\D) \cap H^2(\D)$ and $(-\tfrac{1}{2}\Delta   + V) v \in H^1_0(\D) \cap H^2(\D)$, the approximation properties of the LOD space hence yield
\begin{align} \label{8thOrderError}
\left| \int_\D |u|^2 v -\PLODLtwo(|u|^2)v \dx\right| \le \| |u|^2-\PLODLtwo(|u|^2) \| \, \| v-\PLODLtwo(v) \| \leq CH^8.
\end{align} 
To make this argument rigorous for a particular method requires knowledge of the regularity of $u$ and $v$ (which are typically both related to the exact solution at certain time instances). 
In \cite{Superconv}, it was rigorously demonstrated for a modified Crank--Nicolson method for the time-dependent GPE. 

\section{Fully-discrete approximation of ground states}
\label{section:CompEVP}
In this section we discuss an efficient method for computing ground states in the LOD space. For this, let $\VLOD$ denote the corresponding space, typically constructed based on the bilinear form $a(\cdot,\cdot) = \aV(\cdot,\cdot)$ given by \eqref{choice-LOD-a}. Recalling the setting of section \ref{subsubsection-LOD-ground-states}, we are interested in finding an approximation to the discrete ground state
\begin{align}
\label{LOD-ground-state}
\ugsLOD  =\underset{v \in \mathbb{S} \cap \VLOD}{\mbox{arg\hspace{2pt}min}}\hspace{2pt} \Egs(v).
\end{align}
This can be achieved with a Riemannian gradient method that was originally proposed in \cite{SobolevGradient}. To keep the presentation and the motivation short, we will however introduce the method here as an inverse iteration with damping. Before we formulate a corresponding fully-discrete version of the method in the LOD space, let us first sketch the idea in the idealized analytical setting \eqref{definition-groundstate}, where we seek $\ugs \in H^1_0(\D)$ with
\begin{align*}
\ugs =\underset{v \in \mathbb{S}}{\mbox{arg\hspace{2pt}min}}\hspace{2pt} \Egs(v).
\end{align*}
Recall here that $\mathbb{S}$ is the Riemannian manifold $\mathbb{S}:= \{  v \in H^1_0(\D) \mbox{ and } \| v\| =1 \}$. As the minimization problem can be equivalently written as a nonlinear eigenvalue problem \eqref{definition-groundstate-2} seeking the eigenfunction to the smallest eigenvalue $\lambdags$, a natural choice is the classical inverse iteration method, where the approximation from the previous iteration can be used to linearize the differential operator. For that, let us denote by 
$$
\mathcal{A}_z^{-1} : H^{-1}(\D) \rightarrow H^1_0(\D)
$$ 
the inverse of the linearized Gross--Pitaevskii operator in a linearization point $z \in H^1_0(\D)$. To make the definition of $\mathcal{A}_z^{-1}$ precise, we consider the bilinear form associated with the $z$-linearized Gross--Pitaevskii operator given by
\begin{align}
\label{def-az}
a_{z}(v,w):=  \tfrac{1}{2} \Ltwo{ \nabla v }{ \nabla w } + \Ltwo{ \,(\Vgs\, + \betags |z|^2) \, v }{ w }.
\end{align}
With this, the preimage $\mathcal{A}_z^{-1}F \in H^1_0(\D)$ for some $F \in H^{-1}(\D)$ is given as the solution to the linear elliptic problem
$$
a_{z}( \, \mathcal{A}_z^{-1}F \, , v ) = \langle F , v \rangle \qquad \mbox{for all } v\in H^1_0(\D).
$$
If $\mathcal{I} : H^1_0(\D) \rightarrow H^{-1}(\D)$ denotes the canonical $L^2$-identity given by $\langle \mathcal{I} u , v\rangle := \Ltwo{ u}{ v }$, then the standard inverse iteration would read
$$
u^{(k+1)} := \frac{  \mathcal{A}_{u^{(k)}}^{-1} \mathcal{I} u^{(k)}}{ \|  \mathcal{A}_{u^{(k)}}^{-1} \mathcal{I} u^{(k)} \| }.
$$
However, in our nonlinear setting, this is not sufficient to guarantee global convergence to the ground state. For that we combine the inverse iteration with a damping strategy to obtain the modified iterations as 
\begin{align}
\label{ideal-damped-inverse-iteration}
u^{(k+1)} := \frac{ (1 - \theta_k ) u^{(k)}  \, + \, \theta_k \, \gamma_k \, \mathcal{A}_{u^{(k)}}^{-1} \mathcal{I} u^{(k)} }{ \| (1 - \theta_k ) u^{(k)}  \, + \, \theta_k \, \gamma_k \, \mathcal{A}_{u^{(k)}}^{-1} \mathcal{I} u^{(k)} \| },
\quad 
\mbox{where }
\gamma_k := ( u^{(k)} , \mathcal{A}_{u^{(k)}}^{-1} \mathcal{I} u^{(k)} )^{-1}
\end{align}
and the damping parameter $\theta_k \in (0,2)$ is determined numerically by line search (cf. \cite{SobolevGradient}) such that the energy becomes minimal, i.e.,
\begin{align*}
\theta_k := \underset{\theta>0}{\mbox{arg\hspace{2pt}min}} \hspace{2pt}\Egs \big( u^{(k+1)}(\theta)\big),
\end{align*}
where $u^{(k+1)}(\theta)$ is defined by \eqref{ideal-damped-inverse-iteration} if $\theta_k$ is replaced by the variable parameter $\theta$. The efficient implementation of the line search (which comes at negligible costs) is described in \cite{SobolevGradient}. The appearance of the parameter $\gamma_k$ in addition to the damping parameter $\theta_k$ can be explained from the view point of Riemannian gradient methods, where we again refer to \cite{SobolevGradient} for details. The fundamental advantage of the damped inverse iteration \eqref{ideal-damped-inverse-iteration} is that it converges globally to the ground state for any nonnegative starting value $u^{(0)} \in \mathbb{S}$. Furthermore, the method is strictly energy diminishing (cf. \cite{SobolevGradient}). The guaranteed convergence is a unique feature, which makes the approach very attractive for practical purposes. 

It is straightforward to discretize the damped inverse iterations \eqref{ideal-damped-inverse-iteration} by replacing the analytical functions $u^{(k)}$ by corresponding approximations in the LOD space. However, as discussed in section \ref{subsection:treatment-nonlinearity-with-LOD}, this direct discretization would require the storing of a 4-valence tensor of products of LOD-basis functions, which is computationally not feasible. Hence, we follow the strategy described in section \ref{subsection:treatment-nonlinearity-with-LOD} and replace $|z|^2$ in the definition of $a_z(\cdot,\cdot)$ \eqref{def-az} by $\PLODLtwo(|z|^2)$, where we recall $\PLODLtwo$ as the $L^2$-projection into the LOD space given by \eqref{PLODL2}. 

With this we can now explicitly formulate the fully-discrete method in the LOD space. For that, let $u^{(0)}_{\LOD} \in \VLOD$ be a given initial value and $u^{(k)}_{\LOD} \in \VLOD$ the current iterate. The corresponding linearized bilinear form is given by
\begin{align*}
a_k^{\LOD}(v,w):=  \tfrac{1}{2} \Ltwo{ \nabla v }{ \nabla w } + \Ltwo{ (V + \beta \, \PLODLtwo(|u^{(k)}|^2)) \, v }{ w }.
\end{align*}
To keep the notation short, we denote the discretized inverse (as an approximation of $\mathcal{A}_{u^{(k)}}^{-1}$) by $ \mathcal{D}_{k}^{-1}$, i.e., for $F\in H^{-1}(\D)$, the inverse $\mathcal{D}_{k}^{-1}F \in \VLOD$ is defined as the solution to
$$
a_{k}^{\LOD}(  \mathcal{D}_{k}^{-1}F  , v ) = \langle F , v \rangle \qquad \mbox{for all } v \in \VLOD.
$$
Straightforwardly, we now obtain the next iterate $u^{(k+1)}_{\LOD} \in \VLOD$ by the corresponding discrete damped inverse iteration:
\begin{align}
\label{ideal-damped-inverse-iteration-LOD}
u^{(k+1)}_{\LOD} := \frac{ (1 - \theta_k ) u^{(k)}_{\LOD}  \, + \, \theta_k \, \gamma_k^{\LOD} \, \mathcal{D}_{k}^{-1} \mathcal{I} u^{(k)}_{\LOD} }{ \| (1 - \theta_k ) u^{(k)}_{\LOD}  \, + \, \theta_k \, \gamma_k^{\LOD} \, \mathcal{D}_{k}^{-1} \mathcal{I} u^{(k)}_{\LOD} \| },
\qquad 
\mbox{where }
\gamma_k^{\LOD} := ( u^{(k)}_{\LOD} , \mathcal{D}_{k}^{-1} \mathcal{I} u^{(k)}_{\LOD} )^{-1}
\end{align}
The optimal value of the damping parameter $\theta_k \in (0,2)$ is determined by the property
\begin{align*}
\theta_k := \underset{\theta>0}{\mbox{arg\hspace{2pt}min}} \hspace{2pt}\Egs_{\LOD}\big(u^{(k+1)}_{\LOD}(\theta)\big)
\end{align*} 
for the modified energy
\begin{align} \label{mod_energy}
E_{\LOD}(v) := \int_{\D} \tfrac{1}{2}|\nabla v|^2 + \Vgs |v|^2 + \tfrac{\betags}{2} |v|^2  \PLODLtwo(|v|^2) \dx.
\end{align}
Note that each iteration \eqref{ideal-damped-inverse-iteration-LOD} requires the solving of one elliptic problem involving the computation of $\PLODLtwo(|u^{(k)}|^2)$ and afterwards one more equation to compute $\mathcal{D}_{k}^{-1} \mathcal{I} u^{(k)}_{\LOD}$. Recall here from section \ref{subsection:treatment-nonlinearity-with-LOD} that the computation of $\PLODLtwo(|u^{(k)}|^2)$ is cheap, as it only involves the (nicely conditioned and precomputed) mass matrix in the LOD-space and the precomputed 3-valence tensor. Similarly, all components that are needed to assemble the system matrix  $\mathcal{D}_{k}^{-1} \mathcal{I} u^{(k)}_{\LOD}$ are also precomputed. This makes the iterations in the LOD-space highly efficient. Finally, let us note that the LOD-specific iterations proposed in equation \eqref{ideal-damped-inverse-iteration-LOD} have (in this particular modification) not yet been studied in the literature, neither the convergence in terms of iterations nor the approximation properties of the limit function for $k\rightarrow \infty$. However, our numerical results (cf. section \ref{section:numerical-experiments}) indicate that we can expect the same order of accuracy (cf. section \ref{subsubsection-LOD-ground-states}) as for the ideal LOD ground state $\ugsLOD$ given by \eqref{LOD-ground-state}.

\section{Fully-discrete approximation of the time-dependent problem}
\label{section:high-order-cG-FEM}
In this section we present a fully-discrete method to solve the time-dependent GPE \eqref{model-problem}. For the spatial discretization we use again the generalized finite element space $\VLOD$ based on the Localized Orthogonal Decomposition which was discussed in the previous sections. Hence, we now focus on the time discretization for which we use an energy-conserving time integrator based on continuous Galerkin time stepping. The time stepping method was first introduced and analyzed in \cite{KaM99} for the nonlinear Schr\"odinger equation in two dimensions and generalized in \cite{DoHen22} for the Gross--Pitaevskii equation in up to three dimensions. However, both works assume standard finite element spaces of order $p$ for the space discretization such as the space of $H^1$-conforming Lagrange finite elements. In contrast to that, we now combine this time integrator with the LOD discretization in space. For that purpose, we consider \eqref{model-problem} on a time interval $[0,T]$ and introduce an equidistant time grid $0 = t_0 < t_1 < \dots < t_N = T$ of $n$ time intervals $I_n = (t_n,t_{n+1}]$ with step size $\tau = t_{n+1} - t_n > 0$. Of course, we could also choose a non-equidistant time grid, cf. \cite{DoHen22,KaM99}, but we assume a constant step size $\tau$ for simplicity. Next we choose a polynomial degree $q \in \mathbb{N}$ and define the space
\begin{align*}
	\mathcal{W}_{\tau,q,\LOD}^{H,h,\ell} := \big\{ v : \D \times (0,T] \rightarrow \mathbb{C}\, : \, v_{| \D \times I_n} (x,t) = \sum_{j = 0}^q t^j \, \varphi_{n j} (x),\, \mbox{where }\varphi_{nj} \in \VLOD \big\}.
\end{align*}
Hence, $\mathcal{W}_{\tau,q,\LOD}^{H,h,\ell}$ consists of piecewise polynomials of degree $q$ in time with values in the LOD space $\VLOD$. With this, the numerical solution scheme for the time-dependent GPE, to which we refer as the cG-LOD-method, reads as follows: \\
For a given initial value $u_0 \in H^1_0(\D)$ we seek $u_{\tau}^{\LOD} \in \mathcal{W}_{\tau,q,\LOD}^{H,h,\ell}$ such that for all $n = 0,\dots,N-1$
\begin{align} \label{cGscheme}
\begin{split}
	\int_{I_n} (\ci \partial_t u_{\tau}^{\LOD}, v) - \aV( u_{\tau}^{\LOD}, v) - \beta (\PLODLtwo(|u_{\tau}^{\LOD}|^2) u_{\tau}^{\LOD},v) \, \mathrm{d} t & = 0 \quad \forall v \in \mathcal{W}_{\tau,q-1,\LOD}^{H,h,\ell}, \\
	u_{\tau}^{\LOD,n+} & = u_{\tau}^{\LOD}(t_{n}),
\end{split}
\end{align}
where $u_{\tau}^{\LOD,n+} = \lim_{t \searrow t_n} u_{\tau}^{\LOD}(t)$ and $u_{\tau}^{\LOD,0+} = \PHLOD (u_0)$ and the bilinear form $\aV(\cdot,\cdot)$ is given by \eqref{choice-LOD-a}. We note that according to the previous sections the cG-LOD-method \eqref{cGscheme} proposed here uses the modified nonlinearity $\PLODLtwo(|u|^2) u$ instead of the original cubic nonlinearity $|u|^2 u$. As it was discussed in section \ref{subsection:treatment-nonlinearity-with-LOD}, this decreases the computational effort significantly when assembling the nonlinear terms. Nevertheless, this causes that only the modified energy $E_{\LOD}$ from \eqref{mod_energy} is preserved by the time integration rather than the exact energy $E$. More precise, when testing in \eqref{cGscheme} with $v = \partial_t u_\tau^{\LOD}$ one finds
\begin{align*}
	E_{\LOD}(u_{\tau}^{\LOD}(t_n)) = E_{\LOD}(\PHLOD(u_0)) \qquad \mbox{for all } n = 0,\dots,N.
\end{align*}
For the Crank-Nicolson time integrator it was shown in \cite{Superconv} that this is enough to approximate the exact energy $E(u)$ with an accuracy of $\mathcal{O}(H^6)$ independent on $\tau$. However, a corresponding result for the cG-LOD-method \eqref{cGscheme}, as well as a rigorous error analysis, is still left for future work. \\
Next we describe how the time stepping procedure of the cG-LOD-method \eqref{cGscheme} can be realized in an implementation. For this purpose, the method is reformulated using sufficiently high quadrature rules that integrate the time integrals in \eqref{cGscheme} exactly. We denote by $s_{j}$ and $w_{j}$ for $j = 1,\dots,q$ the nodes and weights of the $q$-stage Gauss-Legendre quadrature rule on the unit interval $[0,1]$ which is exact for polynomials of order $2q-1$. Further, we set formally $s_{0} = 0$ and define the Lagrange basis polynomials $L_{i}$, $i = 1,\dots,q$ of degree $q-1$ associated with $s_{1},\dots,s_{q}$ and the Lagrange basis polynomials $\hat{L}_{i}$, $i = 0,\dots,q$ of degree $q$ associated with $s_{0},\dots,s_{q}$ by
\begin{align*}
	L_{i}(s) = \prod_{j = 1,j\neq i}^q \frac{s - s_{j}}{s_{i} - s_{j}}, \quad \hat{L}_{i}(s) = \prod_{j = 0,j\neq i}^q \frac{s - s_{j}}{s_{i} - s_{j}} \quad \text{for } s \in [0,1].
\end{align*}
The nodes, weights and basis polynomials are easily transformed to $\overline{I_n}$ so that we define further
\begin{alignat*}{2}
	& t_{n,i} := t_n + \tau s_i, \quad && w_{n,i} := \tau w_i, \\
	& L_{n,i}(t_n + s \tau) := L_{i}(s), \quad && \hat{L}_{n,i}(t_n + \tau s) := \hat{L}_{i} (s).
\end{alignat*}
Then $u_\tau^{\LOD}$ is uniquely expressed on $I_n$ via
\begin{align} \label{timesteprepresentation}
	u_\tau^{\LOD} (t) = \sum_{j = 0}^q \hat{L}_{n,j}(t) u^{n,j}_{\tau} \quad \text{for } t \in I_n
\end{align}
with suitable $u_\tau^{n,j} := u^{\LOD}_\tau(t_{n,j}) \in \VLOD$. Testing in \eqref{cGscheme} with $v(t) = L_{n,i}(t) \psi$ for arbitrary $\psi \in \VLOD$, we can rewrite the cG-LOD-method \eqref{cGscheme} as
\begin{align} \label{timestepping1}
	\sum_{j = 0}^q m_{ij} \Ltwo{u_\tau^{n,j}}{\psi} + \ci \, \tau w_{i} \, \aV(u_{\tau}^{n,i}, \psi) + \ci \beta \int_{I_n} L_{n,i}(t) \Ltwo{\PLODLtwo(|u_{\tau}^{\LOD}|^2)u_{\tau}^{\LOD}}{\psi} \mathrm{d} t = 0
\end{align}
for all $\psi \in \VLOD$ and $i = 1,\dots,q$. Here the coefficients $m_{ij}$ are given by
\begin{align*}
	m_{ij} := \int_0^1 \hat{L}'_{j}(s) L_{i}(s) \mathrm{d} s.
\end{align*}
Next we observe that the integrand of the remaining time integral in \eqref{timestepping1} is of degree $4q-1$ in time. Hence the integral is integrated exactly by a $2q$-stage Gauss-Legendre quadrature rule for which we denote by $\tilde{s}_{j}$ and $\tilde{w}_{j}$ for $j = 1,\dots,2q$ the corresponding nodes and weights on the unit interval $[0,1]$. Setting formally $g(u) := P_{L^2}^{\LOD}(|u|^2)u$ we observe, in view of \eqref{timesteprepresentation}, that \eqref{cGscheme} is equivalent to
\begin{align} \label{timestepping2}
		\sum_{j = 0}^q m_{ij} \, \Ltwo{u^{n,j}_{\tau}}{\psi} + \ci \, \tau w_{i}  \, \aV(u^{n,i}_{\tau}, \psi ) + \ci \, \beta \tau \sum_{\nu = 1}^{2q} \tilde{w}_{\nu} \, L_{i}(\tilde{s}_{\nu}) \Big( g \Big( \sum_{j = 0}^q \hat{L}_{j}(\tilde{s}_{\nu}) u^{n,j}_{\tau} \Big), \psi \Big) = 0
\end{align} 
for all $\psi \in \VLOD$ and $i = 1,\dots,q$. Hence in each time step one has to solve for given $u^{\LOD}_\tau(t_n) = u^{n,0}_\tau$ the nonlinear system \eqref{timestepping2} for $u^{n,j}_\tau$, $j = 1,\dots,q$ so that by the continuity condition in \eqref{cGscheme} the solution at the new time instance $t_{n+1}$ is given by
\begin{align*}
	u^{\LOD}_\tau(t_{n+1}) = \sum_{j = 0}^q \hat{L}_{j}(1) u^{n,j}_\tau.
\end{align*}
Note that the the system \eqref{timestepping2} is a fully coupled system for the unknowns $u^{n,j}_\tau$, $j = 1,\dots,q$. However, as proposed in \cite{KaM99}, the system can be decoupled by observing that for $\mathcal{M} := (m_{ij})_{i,j=1,\dots,q}$ and $W := \mathrm{diag}(w_1,\dots,w_q)$ the matrix $\mathcal{M}^{-1} W$ is the coefficient matrix of the $q$-stage Gauss-Legendre Runge-Kutta method. Therefore, it is diagonalizable (cf. \cite{Dekker84}) such that $\Sigma \mathcal{M}^{-1} W \Sigma^{-1} = \Gamma = \mathrm{diag}(\gamma_1,\dots,\gamma_q)$ for some matrix $\Sigma \in \mathbb{R}^{q \times q}$. Defining the new unknowns
\begin{align*}
	U^{n,i} = \sum_{j = 1}^q \Sigma_{ij} u_\tau^{n,j}, \quad i = 1,\dots,q
\end{align*}
in each time step one infers that \eqref{timestepping2} is equivalent to
\begin{align} \label{timesteppingtransform}
	(U^{n,i}, \psi) + \ci \, \tau \gamma_i \, \aV(U^{n,i}, \psi) = A_i(u^{n,0}_{\tau}, \psi) - \ci \, \beta \tau \sum_{\nu = 1}^q B_{i\nu} \Big( g \big( C_{0\nu} u^{n,0}_{\tau} + \sum_{j = 1}^q C_{j\nu} U^{n,j} \big), \psi \Big)
\end{align}
for all $\psi \in \VLOD$ and $i = 1,\dots,q$. Here the coefficients $A_i, B_{i\nu}, C_{0\nu}, C_{i\nu}$ for $i = 1,\dots,q$ and $\nu = 1,\dots, 2q$ are given by
\begin{align*}
	A_i = \sum_{j = 1}^q \Sigma_{ij}, \quad B_{i\nu} = \tilde{w}_\nu \sum_{j = 1}^q (\Sigma \mathcal{M})^{-1}_{ij} L_j(\tilde{s}_\nu), \quad C_{0\nu} = \hat{L}_0(\tilde{s}_\nu), \quad C_{i\nu} = \sum_{j = 1}^q \hat{L}_j(\tilde{s}_\nu) \Sigma^{-1}_{ji}.
\end{align*}
In order to solve the nonlinear system of equations \eqref{timesteppingtransform} in each time step we propose to use a fixed-point iteration. For this, we define $\mathcal{L}_{\LOD}: \VLOD \rightarrow \VLOD$ via $\Ltwo{\mathcal{L}_{\LOD} v}{w} = \aV(v,w)$ so that we can formulate \eqref{timesteppingtransform} as
\begin{align*}
	(I +  \ci \, \tau \gamma_i \mathcal{L}_{\LOD}) U^{n,i} = A_i u^{n,0}_{\tau} - \ci \beta \sum_{\nu = 1}^q B_{i\nu} g \big( C_{0\nu} u^{n,0}_{\tau} + \sum_{j = 1}^q C_{j\nu} U^{n,j} \big), \quad i =1,\dots,q.
\end{align*}
The solution of this equation $U^{n,j}$, $j = 1,\dots,q$, is now approximated via an explicit fixed-point iteration reading for $k \ge 0$ as
\begin{align} \label{timesteppingfixedpoint}
	(I +  \ci \, \tau \gamma_i \mathcal{L}_{\LOD}) U^{n,i}_{k+1} = A_i u^{n,0}_{\tau} - \ci \beta \sum_{\nu = 1}^q B_{i\nu} g \big( C_{0\nu} u^{n,0}_{\tau} + \sum_{j = 1}^q C_{j\nu} U^{n,j}_k \big), \quad i =1,\dots,q.
\end{align}
Here we may choose $U^{n,j}_0 = U^{n-1,j}$ from the previous time step as an initial value for the fixed-point iteration. Note that the operator $(I +  \ci \, \tau \gamma_i\mathcal{L}_{\LOD})$ is independent of the current time step (at least on a equidistant time grid). Hence, the resulting matrix associated with this operator in the LOD basis can be \textbf{LU}-decomposed in a pre-process before the time stepping. 
Using the \textbf{LU}-decomposition the system of linear equations in \eqref{timesteppingfixedpoint} can be solved efficiently.

\section{Numerical experiments}
\label{section:numerical-experiments}
In this section we verify, by means of numerical experiments, the theoretical properties, efficient implementation and promising outlook of the proposed method. The tests are divided into three major parts, namely, first is verified the convergence of the damped inverse iteration in the LOD-space, subsequently the high order convergence of the time integration is verified using a test problem with known analytical solution, and finally the methods are combined to solve a physical setup in $3d$. Additionally, these major parts are subdivided into smaller sections in order to compare different LOD-spaces and make cross comparisons to other spatial discretizations such as spectral and P1-Lagrange elements. Further, we demonstrate how the replacement of the density $|u\LOD|^2$ by its LOD projection $P^\text{ \tiny LOD}_{L^2}(|u\LOD|^2)$ in our methods highly improves their performance.
All CPU times were measured on a laptop with a {\fontfamily{cmtt} 11th Gen Intel® Core™ i7-1165G7 @ 2.80GHz × 8 } processing unit and 64 GB of RAM.

\subsection{Computation of ground states}
\label{section:ground-states}
In this section, we compute the ground state using the damped inverse iteration \eqref{ideal-damped-inverse-iteration-LOD} for a smooth potential and for a discontinuous potential. Throughout we set the tolerance of the damped inverse iteration to $10^{-10}$ in terms of energy.

\subsubsection{Smooth potential}
\begin{figure}[h!]
	\centering
	\begin{minipage}{0.47\textwidth}
		\centering
		\includegraphics[width=1.08\textwidth]{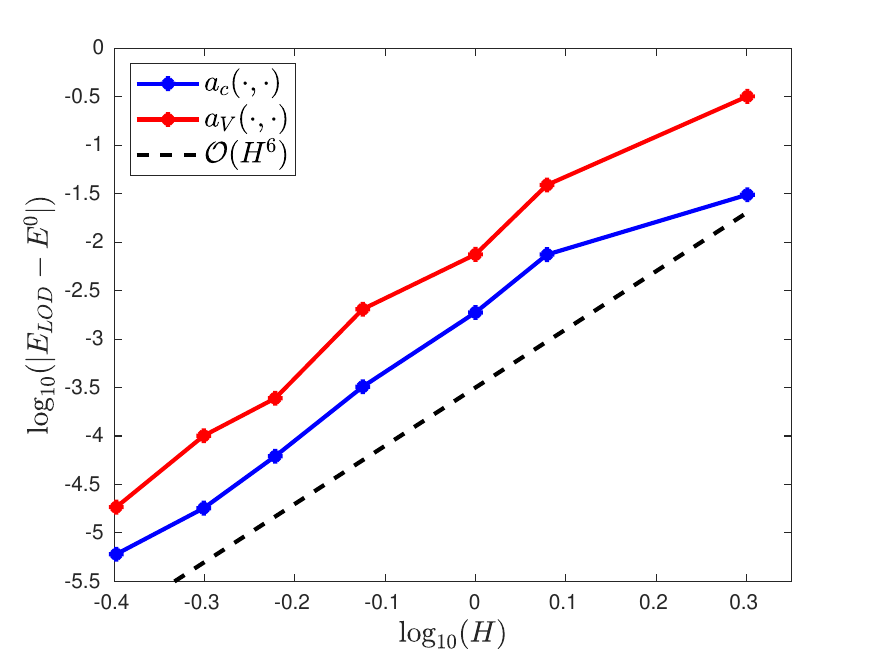} 
		\caption{Convergence plot of values in Table \ref{MinEnergyTable} as measured logarithmically in error versus mesh size $H$}\label{Conv_Plot_2D}
	\end{minipage}\hspace*{\fill}
	\begin{minipage}{0.47\textwidth}
		\centering
		\includegraphics[width=1.08\textwidth]{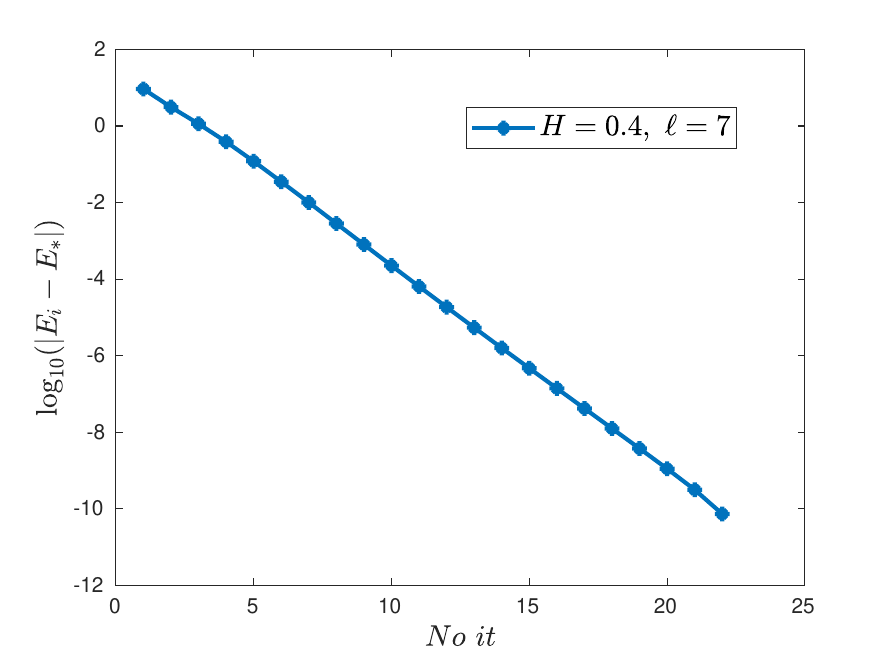} 
		\caption{Convergence to ground state of method \eqref{ideal-damped-inverse-iteration}, initial guess is a normalized Gaussian}\label{Conv_Plot_It}
	\end{minipage}
\end{figure}

\begin{figure}[h!]
	\centering
	\begin{subfigure}{0.50\textwidth}
		\centering
		\includegraphics[width=1.08\textwidth]{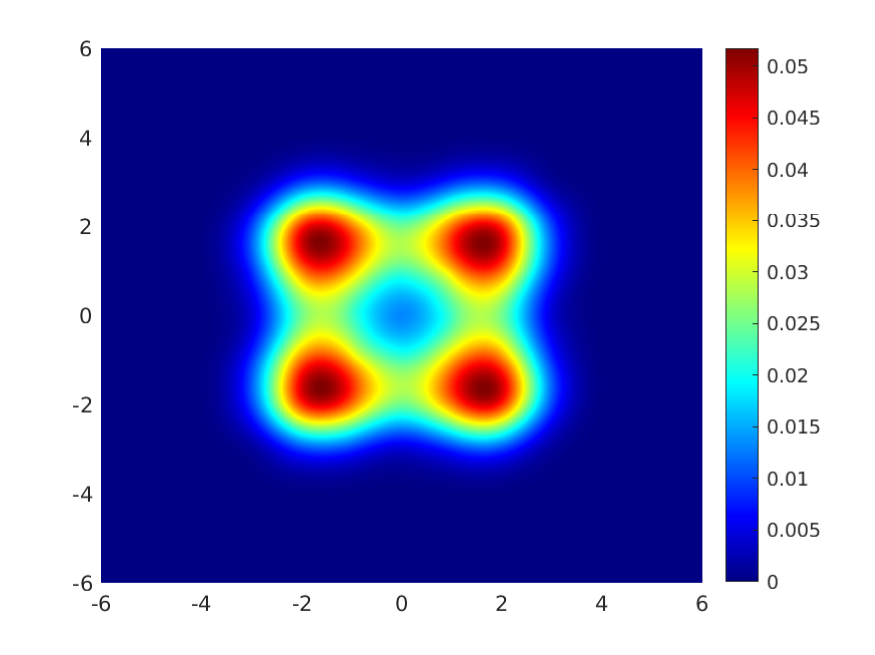} 
		\caption{Canonical choice of inner-product, $\ac$}
	\end{subfigure}\hspace*{\fill}
	\begin{subfigure}{0.50\textwidth}
		\centering
		\includegraphics[width=1.08\textwidth]{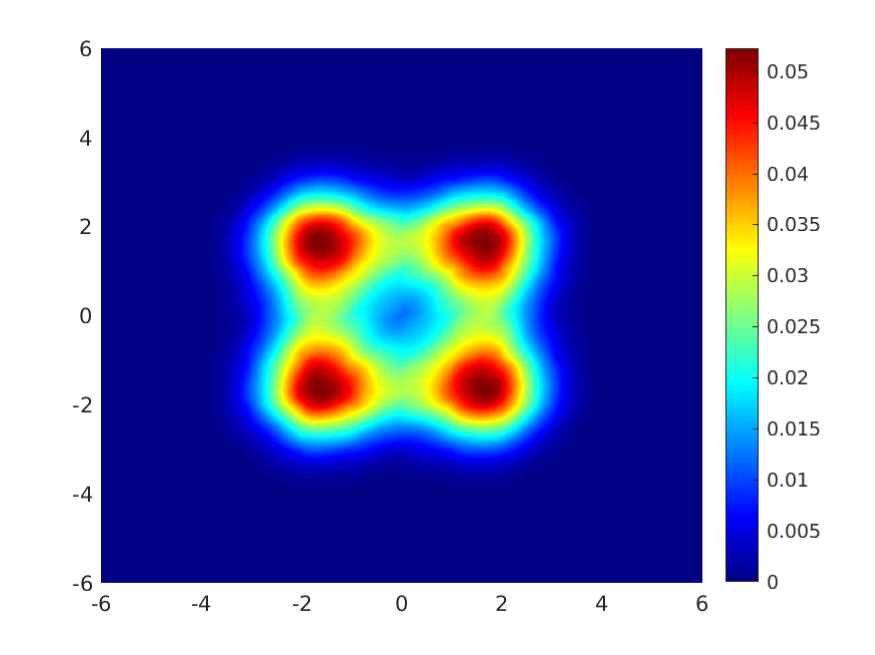} 
		\caption{Problem specific choice of inner-product, $\aV$}
	\end{subfigure}

\caption{Density plot of minimizer, $|\ugsLOD|^2$. At $H=1.0$, $l=3$, the difference in accuracy is more than one order of magnitude.}\label{Comp}
\end{figure}

Our first experiment considers the computation of the minimal energy and the corresponding ground state to a double-well potential in $2d$. More precise, we consider the minimization problem for the ground state
\begin{align} \label{sim-problem1}
	\ugs =\underset{u \in \mathbb{S}}{\mbox{arg\hspace{2pt}min}}\hspace{2pt} \int_\D \frac{1}{2}|\nabla v|^2 + V|v|^2+\frac{\beta}{2}|v|^4 \dx
\end{align}
with $\D = (-6,6)^2$ and where the trapping potential and repulsion parameter are given by
\begin{align} \label{smoothpotential}
	V(x,y) = \frac{1}{2}( x^2+y^2)+4e^{-{x^2}/2}+4e^{-y^2/2} \quad \quad \text{and} \quad \quad \beta = 50.
\end{align}
Recall that $\mathbb{S}$ is the Riemannian manifold $\mathbb{S}:= \{  v \in H^1_0(\D) \mbox{ and } \| v\| =1 \}$. We approximate the ground state solution $\ugs$ of \eqref{sim-problem1}, \eqref{smoothpotential} by a discrete minimizer in a suitable LOD-space given by \eqref{LOD-ground-state} which is then computed with the damped inverse iteration method \eqref{ideal-damped-inverse-iteration-LOD} proposed in section \ref{section:CompEVP}. Besides demonstrating the performance of the damped inverse iteration method in the LOD-space, we also compare the two different choices of the bilinear form for the construction of the LOD-space, i.e.,
\begin{align*}
	\aV(u,v) = \tfrac{1}{2}\Ltwo{\nabla u}{\nabla v} + \Ltwo{V u}{ v} \quad \text{and} \quad \ac(u,v) = \Ltwo{\nabla u}{\nabla v} 
\end{align*}
from \eqref{choice-LOD-a} and \eqref{choice-LOD-a0}. Since $V$ is chosen as a smooth potential we expect that both choices result in the same order of convergence w.r.t. the coarse mesh size, i.e., $\mathcal{O}(H^6)$. After that, we compare the computational cost for our proposed method \eqref{ideal-damped-inverse-iteration-LOD} with the method proposed in \cite{HeMaPe14,HeP21} that does not replace the density $|u\LOD|^2$ by $P^\text{ \tiny LOD}_{L^2}(|u\LOD|^2)$. \\
For the computations the fine mesh size is consistently chosen as $h=H/50$ so as to isolate the convergence in $H$. From a computational point of view, this is too small when $H$ is coarse, since one would ideally choose $h \propto H^3$ so that the error in representing the LOD basis functions on the fine mesh is of the same magnitude as the best approximation error of the solution in the ideal LOD space. However, the choice of the fine mesh $h=H/50$ is chosen to isolate the best approximation error in $H$. \\ 
As a reference solution we choose the result of the damped inverse iteration method \eqref{ideal-damped-inverse-iteration-LOD} with $h = H/80$, $H = 0.3$ and $\ell = 7$.  The energy of the reference solution is computed to be $E^0 = 7.0823112$. Table \ref{MinEnergyTable} shows the obtained results for different values of $H$ (and $\ell$) in terms of the minimal energy together with the estimated errors and the CPU time required in the online phase. The errors are subsequently plotted in Fig. \ref{Conv_Plot_2D}, which clearly demonstrates the $\mathcal{O}(H^6)$ convergence rates of both choices of LOD-space. The canonical construction of the LOD-space with $\ac$ fares better than the problem specific construction incorporating the potential. A striking example of this difference is shown in Fig. \ref{Comp} at which discretization level the difference in accuracy is one order of magnitude. The CPU time in the offline phase is dominated by the computation of the LOD-space and the tensor $\boldsymbol{\omega}$, these values are tabulated in Table \ref{PreComp2d}. For the canonical the LOD-space it suffices, in essence (except for the near $\partial \D$), to compute one basis function which can then be reused and interpolated. We report similar number of iterations for both LOD-spaces and consequently similar CPU times in the online phase as tabulated in the two final columns of Table \ref{MinEnergyTable}.

\begin{table}[h!]
\centering
	\begin{tabular}{| c | c || c  c || c  c || c |c| }
		\multicolumn{2}{c}{} & \multicolumn{2}{c}{$\ac$} & \multicolumn{2}{c}{$\aV$} &\multicolumn{2}{c}{ } \\
		\cline{1-8} 
		$H$ & $\ell$ & $E\LOD$  & $ E\LOD - E^0$&  $E\LOD$  & $ E\LOD - E^0$ & CPU [s] & N$^\text{o}$ it \\
		\hline
		$2.0$ & $1$ & 7.1131522  & 0.0308409&  7.4010560 & 0.3187454 &3 & 18 \\
		$1.2$ & $2$ &7.0897461	  & 0.0074349 &  7.1211954  & 0.0388848  &3 & 18  \\
		$1.0$ & $3$ & 7.0841921  & 0.0018808 &  7.0897783  & 0.0074678 & 4 & 19 \\
		$0.75$ & $4$ & 7.0826330  & 0.0003217 &   7.0843425 & 0.0020320 &4 & 21 \\
		$0.60$ & $5$ &7.0823724  & 0.0000612 &  7.0825553 &  0.0002448  & 7 &21\\
		$0.50$ & $6$ &7.0823286 & 0.0000173 &  7.0824111 &0.0001005 &16&22 \\
		$0.40$ & $7$ & 7.0823166 & 0.0000054  &  7.0823290  & 0.0000184 & 35 &22\\
		$0.30$ & $7$ & 7.0823112 &  $0.000001^*$ \normalfont  &  & &  75& 22 \\ 
		\hline
	\end{tabular}
	\caption{Energy $E\LOD$ from \eqref{mod_energy} of the ground state obtained by \eqref{ideal-damped-inverse-iteration-LOD} for different discretizations using the two different LOD-spaces constructed by $\ac$ and $\aV$. The fine grid discretization was 50 times finer than the coarse, i.e., $h = H/50$, except for the reference energy $E^0$, for which $h=H/80$ and $H = 0.3$. Asterisk marks estimated value. Tabulated CPU time denotes to the total online time.}\label{MinEnergyTable}
\end{table}

\begin{table}[h!]
	\centering
	\begin{tabular}{| c | c || c  c || c |c| }
		\multicolumn{2}{c}{} &  \multicolumn{2}{c}{$\aV$} &\multicolumn{2}{c}{ } \\
		\cline{1-6} 
		$H$ & $\ell$ & $E$  & $E - E\LOD $ & CPU [s] & N$^\text{o}$ it \\
		\hline
		$2.0$ & $1$ &  7.4459532 &   0.044897 &12 & 18 \\
		$1.2$ & $2$ &  7.1305263  & 0.009331  &30 & 18  \\
		$1.0$ & $3$  &  7.0913014  & 0.001523 & 47 & 19 \\
		$0.75$ & $4$&   7.0844227 & 0.000080 &189 & 21 \\
		$0.60$ & $5$&  7.0825835 &  0.000028  & 398 &21\\
		$0.50$ & $6$  &  7.0824042 &-0.000007 &682&22 \\
		$0.40$ & $7$  &  7.0823313 &0.000002 &1651&22 \\
		\hline
	\end{tabular}

	\caption{Energy $E$ from \eqref{GPenergy} for the ground state obtained for different discretizations using the operator-adapted LOD-space without replacement of $|u\LOD|^2$ by $P^\text{ \tiny LOD}_{L^2}(|u\LOD|^2)$, i.e., corresponding to the method analyzed in \cite{HeMaPe14,HeP21}. The values for $E\LOD$ are taken from Table \ref{MinEnergyTable}. The fine grid discretization was 50 times finer than the coarse, i.e., $h = H/50$. Tabulated CPU time denotes to the total online time.}\label{MinEnergyTable2}
\end{table}

Next, we do the same experiment using the operator-adapted LOD-space but without replacement of $|u\LOD|^2$ by $P^\text{ \tiny LOD}_{L^2}(|u\LOD|^2)$, i.e., corresponding to the method analyzed in \cite{HeMaPe14,HeP21}. The results are shown in Table \ref{MinEnergyTable2}. First of all we see that the error caused by replacing the density $|u\LOD|^2$ by $P^\text{ \tiny LOD}_{L^2}(|u\LOD|^2)$ converge to zero in terms of the energy as the coarse mesh size tends to zero. In particular, the error is negligible compared to the error caused by the spatial discretization shown in Table \ref{MinEnergyTable} since it is roughly one magnitude smaller. However, the computational time decreases enormously when including the projected density into the method. Our method proposed in this work is roughly 50 times faster than the method from \cite{HeMaPe14,HeP21}. At the same time this comes with almost no loss w.r.t. its approximation properties, though this still has to be investigated analytically. This indicates a significant boost in efficiency.


\begin{table}[H]
	\centering
	\begin{tabular}{|c| c |c| c|c|c|  }
		\multicolumn{6}{c}{Precomputations} \\ \cline{1-6}
		$H$ & $\ell$ &CPU [s] / basis  & $ \boldsymbol{\omega} $ CPU [s] & mem LOD-matrix $\Phi^{\LOD}$  & mem $\boldsymbol{\omega}$  \\
		\hline
		$2.0$ &  $1$ & 4 & 5 & 11 MB & 11 kB \\
		$1.2$ & $2$ & 8 & 4 & 82 MB & 800 kB  \\
		$1.0$ & $3$ & 13 & 7 & 283 MB&  3.4 MB \\
		$0.75$ &  $4$ & 25 & 18 & 580 MB& 18 MB \\
		$0.6$ & $5$ & 54 &46 & 1.4 GB & 93 MB\\
		$0.5$ & $6$ & 120 &140 & 2.6 GB & 180 MB  \\
		$0.4$ & $7$ & 170 &468 & 5.5 GB & 350 MB \\
		$0.3$ & $7$ & 509 & 729&27.1 GB & 1.2 GB  \\
		\hline
	\end{tabular}
	\caption{Computational times in the offline phase. The basis functions are computed to a relative tolerance of $10^{-7}$. In theory, the number of basis functions necessary to compute with inner product $\ac$ is only of $\mathcal{O}(\ell)$ whereas with $\aV $ it is of $\mathcal{O}(H^{-2})$. Memory usage, all LOD-basis functions are listed as represented on the fine mesh. Similarly, all values of $\boldsymbol{\omega}$ are listed, in spite of many being equal. }\label{PreComp2d}
\end{table}
\subsubsection{Discontinuous potential}
\begin{figure}[H]
	\centering
	\begin{subfigure}{0.45\textwidth}
		\centering
		\includegraphics[width=1.08\textwidth]{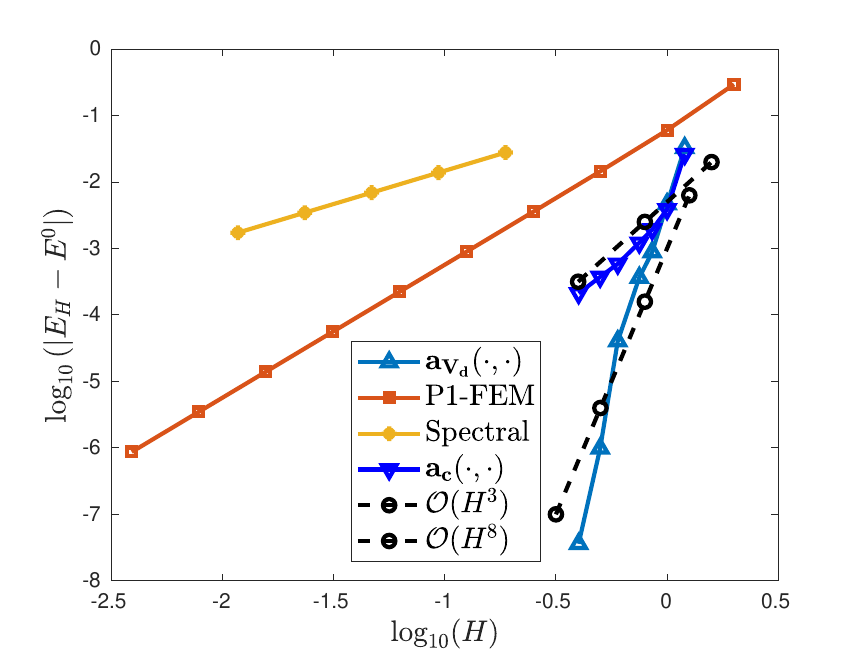} 
		\caption{Logarithmic convergence plot of error versus mesh size $H$. }\label{conv_discont}
	\end{subfigure}\hfil
	\begin{subfigure}{0.45\textwidth}
		\centering
		\includegraphics[width=1.08\textwidth]{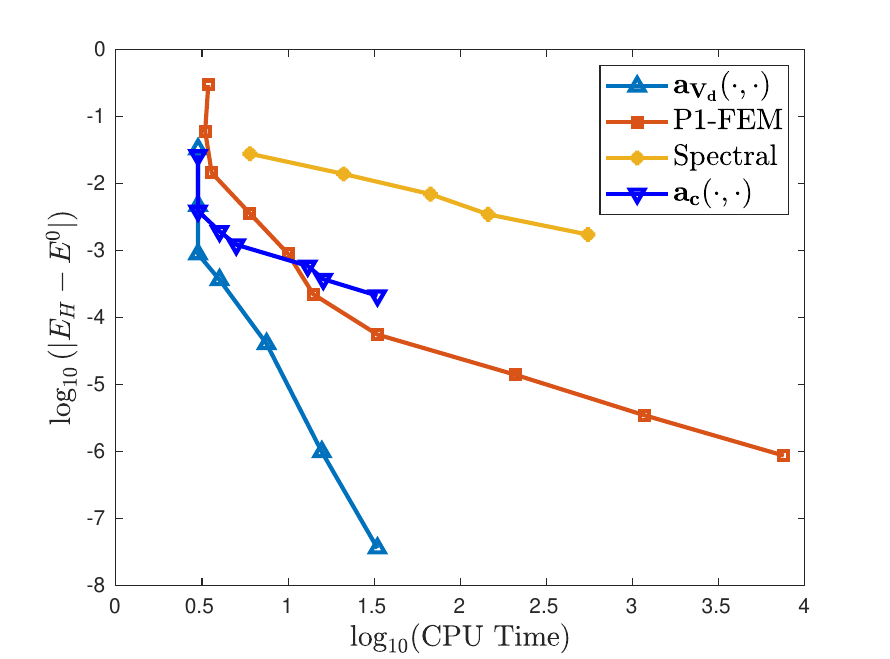} 
		\caption{Logarithmic plot of error vs. CPU time in the online phase.}\label{Online_time_discont}
	\end{subfigure}
	\begin{subfigure}{0.45\textwidth}
		\centering
		\includegraphics[width=1.08\textwidth]{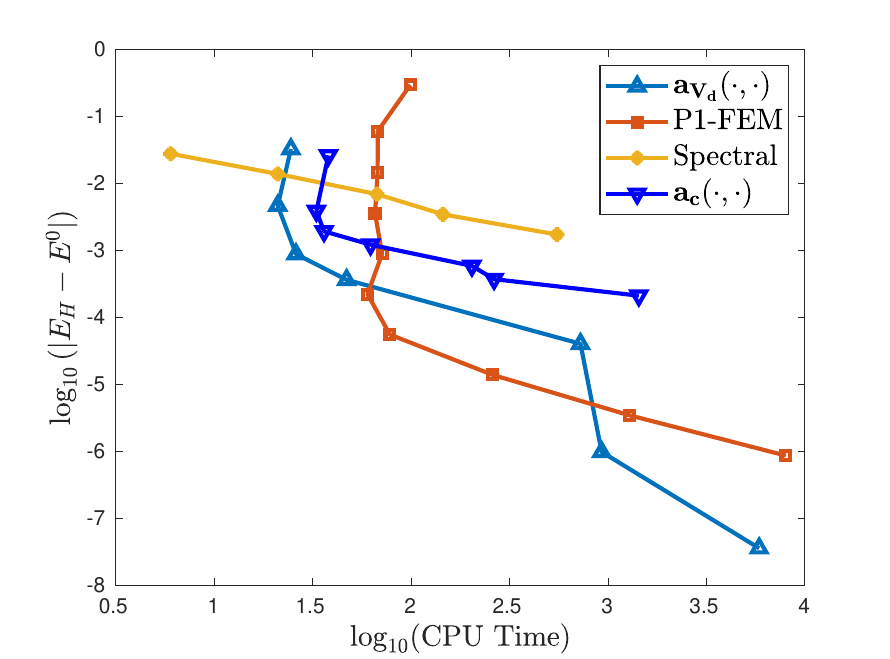} 
		\caption{Log-log plot of error versus total CPU time including computation of basis functions and tensor.}\label{Total_time_discont}
	\end{subfigure}\hfil
	\begin{subfigure}{0.45\textwidth}
		\centering
		\includegraphics[width=1.08\textwidth]{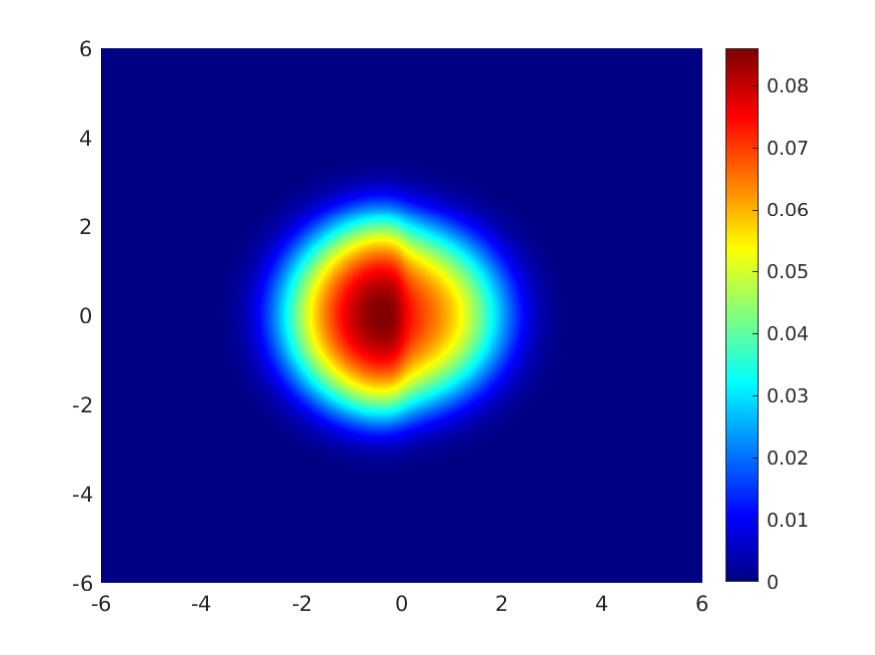} 
		\caption{Density of solution in the LOD-space for $H=0.4$, $Nh=50, \ \ell = 8$.}
	\end{subfigure}
	\caption{Comparison of different methods in terms of convergence and accuracy per CPU times. The blue line with downward pointing triangles in Figures \ref{conv_discont} through \ref{Total_time_discont} show the performance of the canonical LOD. The blue line with upwards pointing triangles shows the performance of incorporating the discontinuity in the LOD-space. The choice of $h$ and $\ell$ is according to the triplets $(H,H/h,\ell)$: (1.2,25,1); (1,25,3);(6/7,25,4);(3/4,25,5);(0.6,40,6);(0.5,40,6);(0.4,50,8). }\label{UGS_discont}
\end{figure}

In the previous test case, nothing was gained by incorporating the potential into the construction of the LOD-space. For smooth problems with analytical solutions the canonical LOD-space is, computationally, a better choice. Furthermore, the smooth potential makes the previous problem well-suited for spectral methods. However, in the case of potentials with only $L^{\infty}$-regularity, such as potentials with a discontinuity, one can at most expect linear convergence w.r.t. energy using a spectral method as the Fourier coefficients of the potential decay as $\mathcal{O}(1/N)$, $N$ being the number of Fourier modes in each dimension. In contrast, the LOD-method achieves 6th order convergence even for $L^2$-potentials, provided the low-regularity potential is incorporated into the construction of the space \cite{HeP21}. To exemplify this we consider the minimization problem
\begin{align} \label{sim-problem2}
	\ugs =\underset{v \in \mathbb{S}}{\mbox{arg\hspace{2pt}min}}\hspace{2pt} \int_\D \frac{1}{2}|\nabla v|^2 + V|v|^2+\frac{\beta}{2}|v|^4 \dx
\end{align}
where again $\D = (-6,6)^2$ and $\beta = 50$ but with a discontinuous trapping potential given by
\begin{align} \label{discontiuouspotential}
	V = V_0 + V_d, \quad V_0(x,y) = \frac{1}{2}( x^2+y^2), \quad V_d(x,y) = \mathbbm{1}_{x \ge 0}(x,y).
\end{align}
Here $\mathbbm{1}_{x \ge 0}$ denotes the indicator function on the set $\{(x,y)\in\R^2 : x>0\}$. Again we approximate the ground state solution $\ugs$ of \eqref{sim-problem2}, \eqref{discontiuouspotential} by a discrete minimizer in the LOD-space given by \eqref{LOD-ground-state} which is then computed with the damped inverse iteration method \eqref{ideal-damped-inverse-iteration-LOD} proposed in section \ref{section:CompEVP}. We construct the LOD-space either with the canonical bilinear form
\begin{align*}
	a_c(u,v) = \Ltwo{\nabla u}{ \nabla u}
\end{align*}
from \eqref{choice-LOD-a} or with
\begin{align*}
	a_{V_d}(u,v) = \tfrac{1}{2}\Ltwo{\nabla u}{ \nabla u} + \Ltwo{\mathbbm{1}_{x \ge 0} u}{v}
\end{align*}
which incorporates the discontinuity of the potential. Due to the symmetry of the discontinuous potential only a handful of basis functions need to be computed, more precisely the $\mathcal{O}(\ell)$ basis functions whose support intersect the line $x = 0$. The rest can be interpolated from these, save for basis functions near the boundary. \\[0.5em]
Next we compare the LOD-method with a standard P1-FEM approach and with the spectral method implemented in the highly successful and widely used package GPELab \cite{gpelab}. A reference value for the minimal energy is computed to 10 digit accuracy using Richardson extrapolation and was found to be $E^0 = 3.341711792$. \\[0.5em]
The convergence rates are illustrated in Fig. \ref{conv_discont}, where downward pointing triangles indicate the canonical LOD-space and upward indicate LOD-space with discontinuity incorporated. The apparent 8th order convergence of the LOD-method is believed to transient and due to the error introduced when replacing the density with its projection onto the LOD-space, i.e., \eqref{8thOrderError}.  In Fig. \ref{Online_time_discont} is plotted the accuracy vs. CPU times in the online phase for the different methods. In this respect the LOD-method outclasses the other methods and we note that in the LOD-space the minimization algorithm converged in around 20 iterations, whereas for the P1 discretization it required around 40. However, due the expensive offline phase, it is not until an accuracy of $10^{-6}$, that something is gained from using the current implementation of the LOD-method, as illustrated in Fig. \ref{Total_time_discont} where the accuracy is plotted versus the total CPU time. Finally, the solution is plotted in Fig. \ref{UGS_discont}.\\
\subsubsection{Harmonic potential}
\begin{figure}[H]
	\begin{subfigure}{0.47\textwidth}
				\centering
		\includegraphics[width=1.08\textwidth]{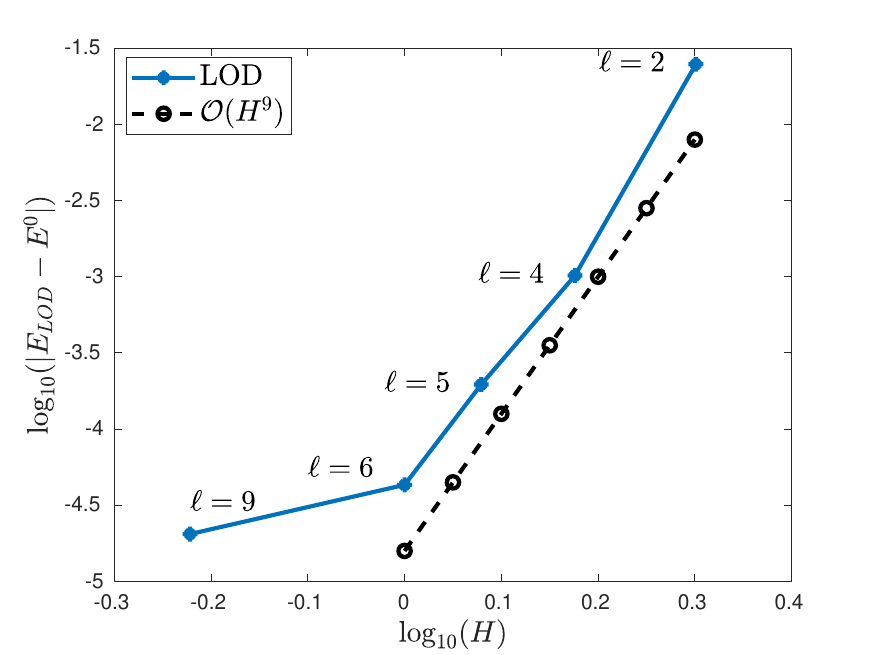} 
		\caption{$d = 2$}\label{Sym_2D}
	\end{subfigure}\hspace*{\fill}
	\begin{subfigure}{0.47\textwidth}
	\centering
	\includegraphics[width=1.08\textwidth]{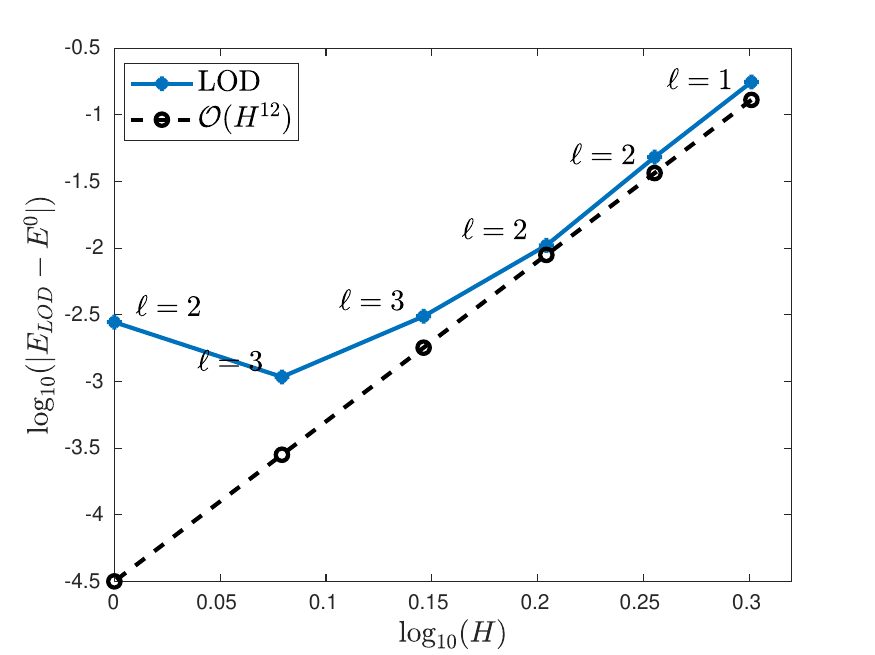} 
	\caption{$d = 3$}\label{Sym_3D}
	\end{subfigure}\hspace*{\fill}
\caption{Log-log plot of error versus coarse mesh size $H$ in the case of harmonic potential.}
\end{figure}

We observe even higher order convergence than the predicted $\mathcal{O}(H^6)$ for purely harmonic potentials. Namely in $2d$ we observe $\mathcal{O}(H^9)$, as illustrated in Fig. \ref{Sym_2D}  and in $3d$, $\mathcal{O}(H^{12})$ is observed, see Fig. \ref{Sym_3D}.  We stress that the observation only holds for approximately global basis functions. In fact, the convergence is so fast that we are unable to conclusively say that this is the asymptotic behavior rather than some incidental initial convergence. This is mostly due to the severe limitations in the offline phase stemming from the representation of the LOD-basis functions as P1-elements. Specifically we consider the following problem:
\begin{align*}
	\ugs =\underset{u \in \mathbb{S}}{\mbox{arg\hspace{2pt}min}}\hspace{2pt} \int_\D \frac{1}{2}|\nabla v|^2 + V|v|^2+\frac{\beta}{2}|v|^4 \dx
\end{align*}
with $V(x) = \tfrac{1}{2}|x|^2$, $x \in \mathbb{R}^d$ and $\beta = 50$.
The domain is $(-6,6)^2$ in $2d$ and $(-5,5)^3$ in $3d$.
Owing to the exponential decay of the ground state, the solution may be well approximated by the solution to the same problem but posed in the full space $\R^3$. Due to the  symmetry, this can be further simplified to a problem in $1d$, namely the following:
\begin{align*}
	\min_{\tilde{u}\in H^1(\R_+)} 2^{d-1} \pi \int_0^\infty\bigg( \frac{1}{2} |\partial_r \tilde{u}|^2 + \frac{1}{2}r^2|\tilde{u}|^2+\frac{\beta}{2}|\tilde{u}|^4\bigg) r^{d-1} dr \quad \text{s.t.} \quad \|\tilde{u}\| = 1,
\end{align*} 
where $d$ denotes the dimension. The corresponding energy can then easily be computed using e.g., a standard FDM such as described in section 3.4 in \cite{BaoNum}. The result, which to 14 digits accuracy is $E^0 = 2.3734292669786$ in $3d$ and $E^0 = 2.896031852200792$  in $2d$, is used as a reference solution.

For brevity, the results are summarized in Fig. \ref{Sym_2D} and in Fig. \ref{Sym_3D}. From these figures it is gathered that in $2d$, $H = 1$, i.e., a grid of $6\times6$, achieves an error of 0.000045. More impressive still is the fact that in $3d$ an error in energy of only $3\cdot 10^{-3}$ is achieved for $H=1.4$.

\subsection{Computations of dynamics}
\label{section:dynamics}
In this section we present numerical experiments for the time-dependent GPE and the combined GPEVP-GPE problem. For the latter, we consider the GPE in three dimensions and simulate the dynamics of BEC in a harmonic trapping potential with the methods presented in this work. Before doing so, we investigate the performance of the LOD discretization combined with the continuous Galerkin time integrator which results in the cG-LOD-method presented in section \ref{section:high-order-cG-FEM}. For the performance test of the cG-LOD-method we choose a one dimensional problem from \cite{Exact2010} that was used for benchmarking in \cite{Superconv}.

\subsubsection{Benchmark problem in $1d$}
\label{section:1d}
The following experiment is only one dimensional, yet it is extremely hard to solve numerically on long time scales. The formal simplicity of the problem also makes it very well suited for benchmarking. The problem is a cubic nonlinear Schr\"odinger equation and takes the form 
\begin{align*}
\ci \partial_t u = - \partial_{xx} u  - 2|u|^2u \qquad \mbox{in } \mathbb{R} \times (0,T],
\end{align*}
with initial value
\begin{align*}
u(x,0)  =  \frac{8(9e^{-4x}+16e^{4x})-32(4e^{-2x}+9e^{2x})}{-128 + 4e^{-6x}+16e^{6x}+81e^{-2x}+64e^{2x} } & .
\end{align*}
As derived in \cite{Exact2010}, the exact solution is given by
\begin{align*}
u(x,t)& =\frac{8e^{4\ci t}(9e^{-4x}+16e^{4x})-32e^{16\ci t}(4e^{-2x}+9e^{2x})}{-128\cos(12t)+4e^{-6x}+16e^{6x}+81e^{-2x}+64e^{2x}}
\end{align*}
and consists of two interacting but standing solitons. What makes it hard to solve numerically is its sensitiveness to energy perturbations. As demonstrated in \cite{Superconv}, any initial error as measured in energy, $E-E_h$, will result in a splitting of the two solitons which travel with a velocity proportional to $\sqrt{|E-E_h|}$ in opposite directions. Consequently there is, on long time scales $T$, a pre-asymptotic regime that is only resolved if the error in energy fulfills $\sqrt{|E-E_h|} \ll T^{-1}$. For example, if the energy converges as $\mathcal{O}(h^2)$ then, refining $h$ by a factor of 2 only decreases the drifting distance between the solitons by a factor of 2.

\begin{figure}[h]
	\begin{subfigure}{0.47\textwidth}
				\centering
		\includegraphics[width=1.08\textwidth]{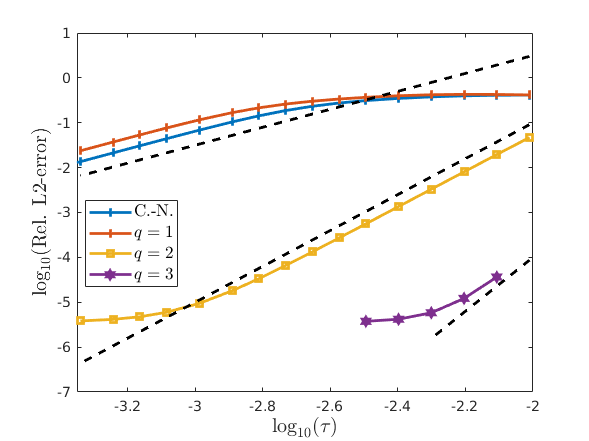} 
		\caption{Convergence plot of rel. $L^2$-error vs. time-step size $\tau$. }\label{q_comp}
	\end{subfigure}\hspace*{\fill}
	\begin{subfigure}{0.47\textwidth}
	\centering
	\includegraphics[width=1.08\textwidth]{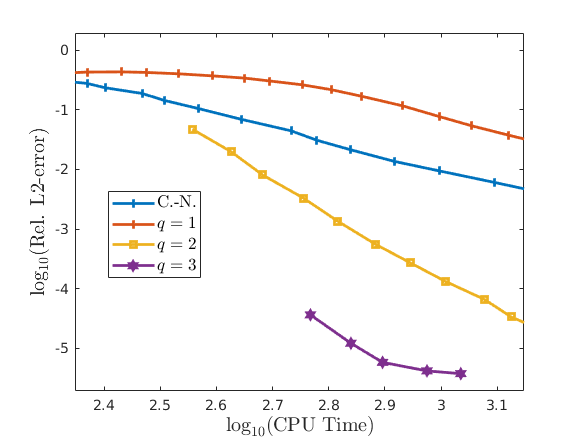} 
	\caption{Convergence plot of rel. $L^2$-error vs. CPU-time. }\label{q_comp_CPU}
	\end{subfigure}\hspace*{\fill}
\caption{Comparison in time of the cG-LOD scheme with $q=1,2,3$ and the Crank-Nicolson scheme. The LOD discretization in space is based on the coarse mesh size $H = 40/2^{11}$ and $\ell = 12$.}
\end{figure}

We consider the fixed spatial discretization $H=40/2^{11}, \ \ell = 12$ on the computational domain $\D = (-20,20)$ with homogenous Dirichlet boundary conditions, motivated by the exponential decay of the solution. First we solve the problem on a small time scale, $T = 10$ and compare the cG-LOD scheme from \eqref{cGscheme} for different choices of $q$ with the standard Crank-Nicolson (CN) time integrator used in \cite{Superconv} which is also based on the LOD discretization in space. The convergence rates for $q=1,2,3$ are plotted in Fig. \ref{q_comp} and we obtain numerically the expected convergence rate in time of $\mathcal{O}(\tau^{2q})$. However, for $q=3$ it is noted that the relative $L^2$-error quickly becomes polluted by the spatial discretization. Nevertheless, as expected the higher order schemes ($q = 2$ and $q = 3$) perform enormously better compared to the CN-scheme in terms of the error and fixed step-size. For $q = 1$ the cG-LOD scheme performs roughly as well as the CN-scheme due to the same order of convergence. It is noted that, the larger $q$ is, the smaller becomes the region of stability of the fixed-point iteration that we proposed in section \ref{section:high-order-cG-FEM} to solve the nonlinear system of equations required in every time step. An example of this is seen in Fig. \ref{q_comp} where a data point is missing for $q=3$ and $\tau = 10/2^{10} \approx 10^{-2}$.

Fig. \ref{q_comp_CPU} shows the relative $L^2$-error versus the CPU time of the cG-LOD scheme for the different choices of $q$ and again compared with the CN-scheme from \cite{Superconv}. When comparing $q = 2$ and $q = 3$ at a fixed computational time of $10^{2.8} \, \mathrm{s} \approx 10 \, \mathrm{min}$ the difference in relative $L^2$-error is of roughly 2 orders of magnitude whereas when comparing $q = 3$ with the CN-scheme we gain more than 3 orders of magnitude in the relative $L^2$-error. A higher order $q$ leads to a gain in CPU-time but imposes a smaller $\tau$. Thus, for a given accuracy one should choose $q$ as large as possible under the constraint that the corresponding time-step size is within the region of stability. For our choice of spatial discretization, selecting $q=3$ is optimal.

\begin{table}[h]
	\centering
	\begin{tabular}{ | l | c | c | c | c |}
	\hline
		\hfil$\tau$ & $\|u-u\LOD\|/\|u\|$ & EOC  &  $\|\nabla(u-u\LOD)\|/\|\nabla u\|$ &  CPU [h]\\ \hline
		$200/22500$ &    0.03841 & - & 0.01619 & 3.2\\
		$200/23000$ &    0.03401 & 5.54&0.01436 &3.3 \\
		$200/23500$&   0.03021 & 5.51 &0.01276 &3.3\\
		$200/24000$&   0.02691 & 5.50 &0.01138 &3.2\\
		
		$200/2^{15}$&   0.00567 & 5.00 &0.00241 &3.9 \\
		$200/(2^{15}\cdot\ 17/16)$&   0.00446 & 3.96&0.00189 &4.1 \\
		$200/(2^{15}\cdot (17/16)^2)$&   0.00362 & 3.44 &0.00154 &4.2 \\
		$200/(2^{15}\cdot (17/16)^3)$&   0.00303 & 2.93 &0.00129 &6.1 \\
		$200/2^{16}$ & 0.00175 &  1.07 & 0.00075 & 6.4 \\
		$200/2^{17}$& 0.00170 & 0.05 &0.00072 & 11.0 \\
		\hline
	\end{tabular}
	\caption{Error table over varying $\tau$ for final time $T=200$, truncation parameter for the LOD is $\ell = 12$ and the coarse mesh size is $H=40/2^{11}$. The estimated order of convergence is computed as EOC$_{i}=\log( E_{i-1}/E_{i})/\log(\tau_{i-1}/\tau_{i})$, where $E_i$ is the $L^2$-error.}\label{Table_Time}
\end{table}

We conclude this one dimensional example with a comparison of previous results from \cite{Superconv}. For that purpose, we now compute the solution on a larger time scale, i.e., we choose $T=200$. Further, we fix $q =3$ and hence solve with the 6th order time-stepping method. The relative errors for different time-step sizes (same as in \cite{Superconv}) are tabulated in Table \ref{Table_Time}. The values demonstrate an initial order of convergence close to 6. However this rate quickly wears off due to the aforementioned spatial pollution. However, compared to the Crank-Nicolson discretization in \cite{Superconv}, which otherwise uses the same approach, we find that a mere 3 hour computation gives the same relative error as the earlier 100 hour computation (which in turn was compared by extrapolation to a hypothetical 4 year computation using a standard $P1$ finite element approach). This significantly improves previous results obtained by Crank-Nicolson time discretizations. Hence, the cG-LOD schemes turns out to be a promising candidate for an efficient solver of the time-dependent GPE in order to simulate dynamics of BECs on long time scales.
 
\subsubsection{Coupled GPEV-GPE problem for $3d$ dynamics}
\label{section:3d}
We conclude this paper with a combined $3$-dimensional GPEV-GPE problem on $\D= (-3,3)^3$, namely:
\begin{align*}
	\begin{alignedat}{2}
		\ci\partial_t u & = -\triangle u +V u+\beta |u|^2u && \qquad \mbox{in } \D \times (0,T], \\
		u  & = 0 && \qquad \mbox{on } \partial \D \times (0,T], \\
		u & = \ugs && \qquad \mbox{on } \D \times \{ t = 0 \},
	\end{alignedat}
	\end{align*}
where $V(x) = \tfrac{1}{2}|x|^2$ and the initial value is given by the ground state w.r.t. to a different trapping potential, i.e.,
\begin{align*}
\ugs = \underset{u \in \mathbb{S}}{\mbox{arg\hspace{2pt}min}}\hspace{2pt} \int_\D \frac{1}{2}|\nabla v|^2+(V+V_d)|v|^2+ \frac{\beta}{2}|v|^4 \dx,
\end{align*}
with $V_d = 2\hspace{-3pt}\mod\hspace{-2pt}(\sum_{i=1}^3 \lfloor x_i\rfloor ,2)$. The coarse mesh size is set to $H=0.2$,  the fine mesh size is set to $H/3$, corresponding to  a grid of $90\times 90\times 90$. The truncation parameter is chosen to be $\ell=2$. With this, the LOD-space is computed w.r.t. the inner product defined by $\tfrac{1}{2}(\nabla u, \nabla v) + (V_du,v)$. Final computational time is set to $T=10$ using a time-step size of $T/2^7=0.078125$. We note that in this setting, we were unable to compute a reference solution on the fine mesh. The problem became only computationally traceable in the LOD space, although we mention that the truncation parameter $\ell = 2$ is chosen to be small to keep the computations feasible. However, in Section \ref{outlook} we discuss further potential improvements of the methods that may serve in future for even better approximations in 3d.

\begin{figure}[H]
	\centering
	\begin{subfigure}{0.49\textwidth}
		\centering
		\includegraphics[width=0.8\textwidth]{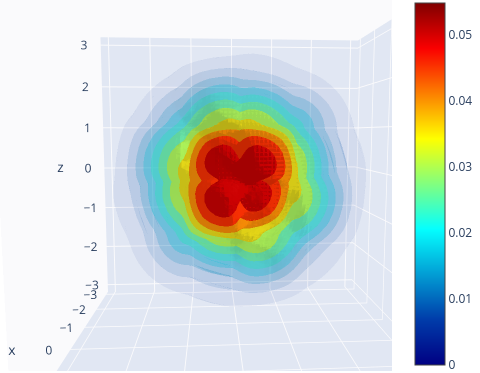} 
		\caption{$t=0$}\label{GS}
	\end{subfigure}
	\begin{subfigure}{0.49\textwidth}
		\centering
		\includegraphics[width=0.8\textwidth]{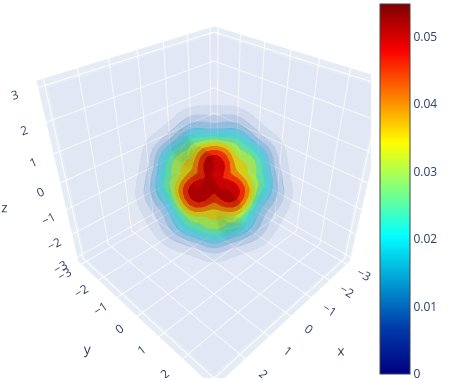} 
		\caption{$t=0$}\label{GS2}
	\end{subfigure}
	\begin{subfigure}{0.49\textwidth}
		\centering
		\includegraphics[width=0.8\textwidth]{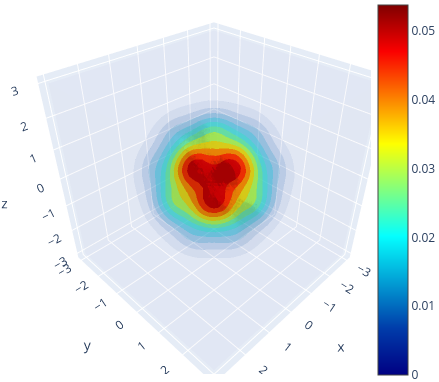} \caption{$t=0.273$}\label{Density0}
	\end{subfigure}
	\begin{subfigure}{0.49\textwidth}
		\centering
		\includegraphics[width=0.8\textwidth]{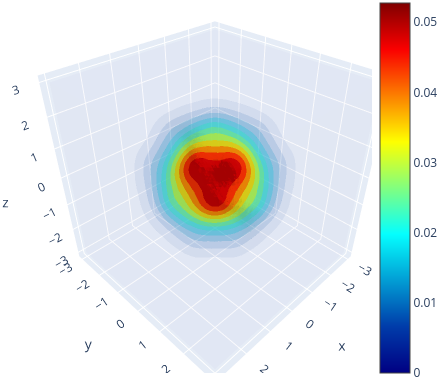} \caption{$t=9.180$}\label{Density1}
	\end{subfigure} 
	\caption{Isosurfaces of the density for model problem in Section \ref{section:3d} at different time instances.}
\end{figure}

We compute the ground state energy to be 3.354636 and the energy for $t>0$ to be 2.425574. Some isosurfaces of the density of the ground state are plotted in Fig. \ref{GS} in which a tetrahedral structure clearly appears. Solving the time-dependent equation, we find that the solution exhibits a periodic behaviour whereby the tetrahedral formation quite abruptly changes orientation. These changes are interlaced with longer periods without such clear structure. This is illustrated in Fig. \ref{GS} through Fig. \ref{Density1}. Comparing for example Fig. \ref{GS} to Fig. \ref{Density1} we see that this abrupt change takes place in around 0.3 units of time. As illustrated in  Fig.\ref{Density1} the density at $t=9.180$ exhibits the same structure as roughly 9 units of time before.

Finally, we give some notes on the CPU times and memory usage: The minimizer in the LOD-space was computed in 690 seconds and required 14 iterations. The solution at final time $T=10$ was computed in about $14$ hours, corresponding to one time step every 400 seconds. As for the offline phase, the LOD-space was computed in about half an hour and the computation of the tensor $\boldsymbol{\omega}$ required 3h. We stress that this was without using any of the periodicity of the discontinuous potential, i.e., many similar basis functions were recomputed and likewise for the entries of $\boldsymbol{\omega}$. The tensor $\boldsymbol{\omega}$ required 22 GB of memory. The representation of the LOD-basis functions on the fine mesh required a matrix of 1 GB.

\section{Outlook} \label{outlook}
We believe the above numerical examples to clearly illustrate the strong capabilities and promising outlook of the proposed method of combining a two level discretization based on the method of LOD with a high order continuous Galerkin time integrator. It achieves high order accuracy and computational efficiency in settings of both low and high regularity. As a final remark, we mention the possibly tremendous improvement of the method by employing a novel construction of the same LOD-space introduced by M. Hauck and D. Peterseim in \cite{HaP23}. This novel construction achieves a much more  localized basis. Here we recall that the computational complexity of our method depends heavily on the localization parameter, $\ell$, of the LOD-basis functions. For $\ell$ = 0, the resulting operators contain neighbor-to-neighbor interactions, e.g., tridiagonal matrices in $1d$. For example, for our benchmark problem in $1d$, the  localization parameter was $\ell= 12$. As shown by Hauck and Peterseim, in $1d$, the very same LOD-space is spanned by a set of basis functions for which $\ell = 1$. In addition, they show that in $2d$ and $3d$ the decay of the LOD-basis functions can in fact be made to be super exponential, more precisely $\exp(-c\ell^{2})$ in $2d$ and $\exp(-c\ell^{1.5})$ in $3d$, as compared to the purely exponential decay of the construction we use. This may well be game changing and shall naturally be considered in future work.

\section*{Dedication}
This article is dedicated in memory of Professor Assyr Abdulle whose research did not only have huge impact on the field of numerical analysis, but who was  also an inspiration to his students and colleagues. Among various other topics, he was especially fascinated by multiscale problems and their numerical treatment, where he made many significant contributions that formed and steered the direction of the research field over the last years \cite{Abd05,Abd09,AbArPa19,AbdBai12,Abd03,AEEV12,AbdPou16,AbdSchwab04}.
The second author (PH) had the chance to join his working group at the EPFL Lausanne in 2014 as a postdoctoral researcher and to enjoy his highly motivating and educative mentorship. During that time, the shared interest in numerical multiscale methods lead to three joined works  \cite{AbH15,AbH17,AbH17b} on the Localized Orthogonal Decomposition, the same methodology used here for the Gross--Pitaevskii equation, showing how his wide research on multiscale methods also had an influence on this article.

\appendix

\section{Notes on the assembly of the 3-valence tensor $\boldsymbol{\omega}$}
\label{section:implementation-notes}

The 3-valence tensor introduced in section \ref{subsection:treatment-nonlinearity-with-LOD}, is extremely sparse. Therefore, some care is required when implementing it. In this section, we briefly describe an efficient digital representation of $\boldsymbol{\omega}$ and implementation of pertaining algorithms. First we recall the definition of the 3-valence tensor from \eqref{def_omega}:
\begin{align*}
\boldsymbol{\omega}_{ijk} := \Ltwo{ \varphi_i^{\LOD}\varphi_j^{\LOD}}{ \varphi_k^{\LOD}}, \quad 1 \le i,j,k \le N_H.
\end{align*}

\subsection{Compressed sparse row representation of $\boldsymbol{\omega}$}
The sparse representation of $\boldsymbol{\omega}$ is similar to the CSR format of sparse matrices. Thus it consists of four vectors: $\Iptr $, $\text{J}_{\boldsymbol{\omega}}$, $\text{K}_{\boldsymbol{\omega}}$ and $\text{V}_{\boldsymbol{\omega}} $. The information about the non-zero elements of $\boldsymbol{\omega}$ with first index $i$ is found at indices $\Iptr\text{[i]}$ to $\Iptr[\text{i}+1]-1$, in the remaining vectors. That is, defining \sub$_i = \Iptr[i]:\Iptr[i+1]-1$, the non-zero values with first index i have second index $\J[\text{\sub}_i]$, third index $\K[\text{\sub}_i]$ and take values V$_{\boldsymbol{\omega}} [\text{\sub}_i]$. Moreover, the list $\J[\text{\sub}_i]$ is sorted. Searching the list $\J[\text{\sub}_i]$ for values equal to $j$, one is left with a subset of indices indicating elements of $\boldsymbol{\omega}$ that have first index $i$ and second index $j$, let us call this list \sub$_{ij}$. Again, these elements have third index according to the list $\K[\text{\sub}_{ij}]$, which is, in turn, sorted. This allows for the quick access to, and manipulation of, entries by means of successive binary searches. For example, consider the following snippet of code which returns the value of $\boldsymbol{\omega}_{ijk}$, 
\begin{algorithm}[H]
	\caption{Accessing values of $\boldsymbol{\omega}$}
\begin{minted}[mathescape]{julia}
function get_val(W,i,j,k)
	sub_i = W.Iptr[i]:W.Iptr[i+1]-1;
	sub_ij = (sub_i[1]-1).+searchsorted(W.J[sub_i],j) #binary search
	sub_ijk = (sub_ij[1]-1]).+searchsorted(W.K[sub_ij],k) #binary search
	return W.V[sub_ijk]
end
\end{minted}
\end{algorithm}
\subsection{Preallocating $\boldsymbol{\omega}$}
How large should the vectors $\Iptr $, $\text{J}_{\boldsymbol{\omega}}$, $\text{K}_{\boldsymbol{\omega}}$ and $\text{V}_{\boldsymbol{\omega}} $ be? Here we describe an efficient way to determin their size. First, note that the mass matrix in the LOD-space,
$M\LOD$, tells us which pairs of basis functions interact. Here, \quotes{interact} means their $L^2$ inner product is non-zero. Now consider basis function $i$, which interacts pairwise with other basis functions according to the list of non-zero elements of the  column $i$ of $M\LOD$. Let us call this list $s$ and its size $n_s$. It is now straightforward to find all basis functions such that $\int \varphi_i\varphi_j\varphi_k dx \not= 0$, by looping over the list $s$. 
\begin{algorithm}[H]
	\caption{Preallocating $\boldsymbol{\omega}$}
\begin{minted}[mathescape]{julia}
INPUT: N #Dimension of LOD-space
       M_LOD #Mass matrix in LOD-space
N_omega_count = 0
W.Iptr = zeros(N+1)
W.Iptr[1] = 1;
for i = 1:N #N dimension of V_LOD
    s = M_LOD[:,i].nzind;
    n_s = length(s)
    start = searchsorted(s,i)[1]
    for j = start:n_s
        for k = j:n_s
            if(abs(M_LOD[s[j],s[k]])>10^-16) N_omega_count += 1; end
        end
    end
    W.Iptr[i+1] = N_omega_count+1;
end
\end{minted}
\end{algorithm}
The above code snippet can then be modified to fill in the values of  $\text{J}_{\boldsymbol{\omega}}$ and $\text{K}_{\boldsymbol{\omega}}$ by setting {\fontfamily{lmtt}\selectfont J[N\_omega\_count]=s[j]} and likewise for $\text{K}_{\boldsymbol{\omega}}$, in the inner loop.

\subsection{Computing elements of $\boldsymbol{\omega}$}

In the spirit of FEM, the values of $\boldsymbol{\omega}$ are computed in a loop over the coarse simplices. For each simplex, quadrature points of a sufficiently high order rule are computed. Each such point is located in a fine simplex. Using this fine simplex' nodes it is straightforward to see, via the LOD-matrix, which basis functions have support in it. These are then listed in {\fontfamily{lmtt}\selectfont active\_basis} and their respective values, in the quadrature point multiplied by the appropriate quadrature weight, are stored in {\fontfamily{lmtt}\selectfont quad\_vals}. These values are then added to $\boldsymbol{\omega}$ in Algorithm \ref{alg:cap}.

\begin{algorithm}[h]
	\caption{Computing elements of $\boldsymbol{\omega}$ }\label{alg:cap} 
	The following is done for each quadrature point on each coarse simplex.
	\begin{minted}[mathescape]{julia}
INPUT: active_basis #list of LOD-basis functions with support on simplex  
       n_basis_loc  #length of "active_basis"
       quad_vals    #eval. of "active_basis" in quadrature point on simplex
       W            # CSR lists of w
       w_H          #quadrature weight

function Add_to_w(quad_vals,w_H,n_basis_loc,active_basis,W) 
    for i = 1:n_basis_loc
        i_glob = active_basis[i];
        sub_i = W.Iptr[i_glob]:( W.Iptr[i_glob+1]-1);
			
        #Evaluate first basis function in quadrature points
        W1 = quad_vals[i]*w_H;
        for j = i:n_basis_loc
            j_glob = active_basis[j]
            sub_ij= (sub_i[1]-1).+searchsorted(W.J[sub_i],j_glob)
			
            #Evaluate second basis function in quadrature points
            W2 = W1*quad_vals[j]
            for k = j:n_basis_loc
                k_glob = active_basis[k];
                sub_ijk = (sub_ij[1]-1).+searchsorted(W.K[sub_ij],k_glob)[1]
                W3 = W2*quad_vals[k]
                W.V[sub_ijk]+=W3;
            end
        end
    end
end
\end{minted}
\end{algorithm}
Regarding the choice of quadrature we remark that in the case of the canonical LOD-space in $1d$, a 9th order quadrature rule computes the elements of $\boldsymbol{\omega}$ exactly. To see this we let $\triangle^{-1} : L^2((0,1))\mapsto H^1_0((0,1))$ denote the solution operator such that $a_c(\triangle^{-1}f,v)=(f,v) \ \forall v \in H^1_0((0,1))$ and let $W=\ker(P_H)$. Consequently, for $v_H\in V_H$, the following holds true $a_c(\triangle^{-1}v_H,w)=( v_H,w) =0 $. 
	The set  $\{v\in H^1_0((0,1)): v =\triangle^{-1}v_H, \text{ for some } v_H\in V_H\}$ is thus the $a$-orthogonal complement of the kernel of $P_H$, in other words it is the LOD-space as defined in equation \eqref{LODdef}. Lastly, since $v''\LOD$ is piecewise linear, any function $v\LOD\in V\LOD$ must be a piecewise polynomial of order 3.  This observation is due to Hauck and Peterseim \cite{HaP23}. However, we stress that, in higher dimensions it no longer holds that the canonical basis is a piecewise polynomial of order 3. 
 Inspired by this, we therefore employ, in $2d$, the 19 point quadrature rule of 9th order as described in \cite{dunavant1985high}, so as to isolate the 6th order convergence in $H$. This quadrature rule achieves the optimal number of points required for 9th order convergence and has positive weights in the interior of the simplex.
Since we are not aware of any adequate quadrature rule of 9th degree in $3d$, we instead employ a 8th degree tetrahedral quadrature in the $3d$ setting. In this case, the minimal number of quadrature points is known to be 40 (cf. \cite{Solin}) which however has not been achieved yet, to the best of our knowledge. We therefore employ the quadrature rule described in \cite{keast1986moderate} requiring 45 quadrature points. Here we remark an error in paper \cite{keast1986moderate}, namely a minus sign is missing in Table 1, column $W$, row $N=8, M=1$ and the entry should read $-0.393270066412926145d-1$.

\end{document}